\documentclass[a4,12pt]{article}
\advance\topmargin -0.8cm \advance\oddsidemargin -0.6cm
\advance\evensidemargin -0.6cm

\usepackage{amsfonts}
\usepackage{amsmath,ae,subfigure}
\usepackage[T1]{fontenc}
\usepackage{graphicx}
\usepackage{enumerate}
\usepackage{float}

\def \R{\mathbb{R}}
\def \N{\mathbb{N}}
\def \supp{\mbox{supp}}
\def \C{\mathbb{C}}
\def\v{\varepsilon}

\def \un{\hbox{1\!\!I}}

\newtheorem{thm}{Theorem}[section]
\newtheorem{lem}{Lemma}[section]

\newtheorem{defi}{Definition}[section]

\begin{document}

\title{Estimation of a semiparametric  contaminated\\ regression model}
\author{Pierre Vandekerkhove }

\date{\small \it LAMA - CNRS UMR 8050, \\
Universit\'e Paris-Est \\
Pierre.Vandekerkhove@univ-mlv.fr}
\maketitle
\begin{abstract}
\baselineskip 16pt
We consider in this paper a contamined regression model  where the  distribution of the contaminating component
is known when  the Euclidean parameters of the regression model, the noise distribution, the  contamination ratio
and the distribution of the design data are unknown. Our model is said to be semiparametric in the sense that  the probability
density function (pdf) of the noise involved in the regression model is not supposed to belong to a parametric density family. 
When the pdf's of the noise and the contaminating  phenomenon are supposed to be symmetric about zero, we propose
an estimator of the various (Euclidean and functionnal) parameters  of the model, and prove under mild conditions
its convergence. We prove in particular that, under technical conditions all satisfied in the Gaussian case,  the Euclidean part of the model is estimated at the rate $o_{a.s}(n^{-1/4+\gamma})$,  $\gamma>0$.
We recall that, as it is pointed out in  Bordes and Vandekerkhove \cite{BV10},  this result cannot be ignored to go further in the  asymptotic theory for this class
of  models.  Finally the implementation and numerical performances of our method are discussed on several toy examples.
\end{abstract}


\section{Introduction}
Let $(U_i)_{i\geq 1}$  be a sequence of independent and identically distributed (iid)  random variables 
according to a {\it{Bernoulli}} distribution with parameter $p\in(0,1)$. We  consider an  iid sample $(Z_1, \dots, Z_n)$ where for all
$i=1,\dots, n$, $Z_i=(X_i,Y_i)$  is a bivariate random variable  defined, relative to $U_i$,  as follows
\begin{eqnarray}\label{model_array}
\left\{
\begin{array}{c}
 Y_i=a_0+b_0X_i+\varepsilon_i^{[0]},\quad\quad \mbox{if }\quad  U_i=0,\\
Y_i=a_1+b_1 X_i+\varepsilon_i^{[1]},\quad\quad \mbox{if }\quad  U_i=1,\\
\end{array}
\right.
\end{eqnarray}
where the design sequence $(X_i)_{i\geq 1}$, respectively  the errors $(\varepsilon_i^{[j]})_{i\geq 1}$, $j=0,1,$  is a sequence of iid random variables  with cumulative distribution function (cdf) $H$, resp. $F_j$, and probability density function (pdf), $h$, resp. $f_j$,  $ j=0,1$. We suppose in addition that the design sequence is independent  from the errors.  
This model, called  the 2-{\it mixture of regression model},   belongs to the wide class of mixture of regression models which has been studied in \cite{ZZ04};  see also 
\cite{YLN03} in a LOS (length of stay)  medical problem,  \cite{BJ10} for prediction, or  \cite{YH10} in a nonparametric modelling context. Recently Martin-Magniette {\it et  al.} \cite{MMBR08} introduced this model in microarray analysis for the study of the  two color ChIP-chip experiment. Briefly, the Chromatin immunoprecipitation (ChIP) is a well established procedure to investigate proteins associated with DNA. ChIP on chip involves analysis of DNA recovered from ChIP experiments by hybridization to miccroarray. In a two color ChIP-chip experiment, two samples  are compared: DNA fragments crosslinked to a protein of interest (IP) and genomic
DNA (input). The goal is then to identify actual binding targets of the IP, \textit{i.e.} probes whose IP  signal  is significantly larger than the input signal. In the model proposed by Martin-Magniette {\it et al.} \cite{MMBR08} the components of the random vector $ Z_i=(X_i,Y_i)$, see model (\ref{model_array}), corresponds respectively  to the log-input and log-IP intensities of probe $i$ when the (unknown) status of the probe is characterized through a label $U_i$ which is 1 if the probe is {\it enriched} and 0 if it is {\it standard} (not enriched).  Note also that the assumption made by these authors on the   error sequences $(\varepsilon_i^{[j]})_{i\geq 1}$, $j=0,1$,  is that  $\varepsilon_i^{[j]}=\varepsilon_i$ for all $(j,i)\in \left\{0,1\right\}\times \N^*$ where $\varepsilon_i$  is a Gaussian random variable with mean 0 and variance $\sigma^2$ (homoscedaticity with respect to the probe status $U_i$). 

In this work, we propose to weaken this last assumption while completely specifying  the regression model under the probe standard  condition  (the parameter $\theta^{[0]}:=(a_0,b_0)\in \R^2$ and $f_0$ are   supposed to be entirely  known).  Note that this kind  of  assumption arises  naturally  in microarray analysis, see model (\ref{BDV06model})   and references \cite{BH95}, \cite{E07}, or \cite{BDV06} p. 744 formula (22),  where analytic expression of $f_0$, characterizing probe expressivity levels under a certain standard condition, is assumed to be available (generally derived from training data and probabilistic computations).
In particular we will suppose that, in model (\ref{model_array}),  the distribution of the $\varepsilon_i^{[1 ]}$ is seen as a nuisance parameter (it is no longer supposed to belong to a parametric distribution family), turning model (\ref{model_array}) into a purely semiparametric model.  
Note that when $\theta^{[0]}$ is known the  observations  $Y_i$, for $i=1,\dots,n$,  can be centered according to $Y_i:=Y_i-(a_0+b_0X_i)$  which implies a  simplification of model (\ref{model_array}), since we then have
\begin{eqnarray}\label{model_arraybis}
\left\{
\begin{array}{lll}
Y_i&=\varepsilon_i^{[0]}, &\quad \mbox{if }\quad  U_i=0,\\
Y_i&=\alpha+\beta X_i+\varepsilon_i^{[1]}, &\quad \mbox{if }\quad  U_i=1,\\
\end{array}
\right.
\end{eqnarray}
where $\alpha:=a_1-a_0$ and $\beta:=b_1-b_0$.   We suppose in model (\ref{model_arraybis}),  which is from now on our model of interest, that  the  $Z_i=(X_i,Y_i)$'s  distribution admits a pdf  with respect to the Lebesgue measure on $\R^2$
 defined by:
\begin{eqnarray}\label{model}
g(x,y)&=&h(x)g_{Y|X=x}(y)\nonumber\\
&=&h(x)[pf(y-(\alpha+\beta x))+(1-p) f_0(y) ],\quad \quad 
(x,y)\in \R^{2},
\end{eqnarray}
where $f$ denotes the unknown pdf of the $\varepsilon_i^{[1]}$, $f_0$ the known pdf of the $\varepsilon_i^{[0]}$, $h$ the unknown pdf of the $X_i$, $f$ and $f_0$ being supposed to belong to the class of even densities. We will  finally denote  by $\vartheta:=(p,\alpha,\beta)\in (0,1)\times \R^2$  the unknown Euclidean  parameter  of model (\ref{model}). Model (\ref{model_arraybis}) corresponds
exactly to a contaminated  version of the semiparametric additive regression model studied in \cite{Cuz92a}, \cite{Cuz92b} and more recently in  \cite{YMP10}.
On the other hand model (\ref{model}) extends for the first time  to the bivariate case,   the class of semiparametric mixture models introduced  by Hall and Zhou \cite{HZ03} for $\R^s$-valued observations with  $s\geq 3$, and studied later in the univariate case,  through two specific models:
\begin{eqnarray}\label{BMV06model}
g(y)=p f(y-\mu_1)+(1-p)f(y-\mu_2), \quad \quad y\in \R,
\end{eqnarray}
where $(p,\mu_1,\mu_2)\in (0,1/2)\times \R^2$ and $f$, supposed to be even, are unknown, see  \cite{BMV06},  \cite{HWH07},  \cite{MS10}, and 
\begin{eqnarray}\label{BDV06model}
g(y)=pf(y)+(1-p)f_0(y-\mu_2), \quad \quad y\in \R,
\end{eqnarray}
where $(p,\mu)\in(0,1)\times \R$ and  $f$  are unknown,  $f_0$ is known, and the pdfs $f$ and $f_0$ are supposed to be even,  see \cite{BDV06},  \cite{BV10}.\\

The paper is organized as follows. In Section 2  we  present  an M-estimating method, inspired by  \cite{BMV06}, \cite{BDV06} and \cite{BV10},
that allows us to estimate the Euclidean and the functional parameters of model  (\ref{model_arraybis});  in Section 3 we address the semiparametric identifiability problem associated to expression (\ref{model})  and establish rates of convergence of our estimators; in Section 4 we discuss    the performance of our method on simulated  examples and focus our attention on 
the optimization problems encountered during its implementation. When technical results are relegated to the appendix, which corresponds to Section 5.

\section{Estimating method}
In the spirit of \cite{BMV06}, \cite{BDV06} and \cite{BV10}, we will suppose that $f$ and $f_0$
are both  pdfs symmetric  about zero (recall  that only  $f_0$ is assumed  known). To avoid trivial situations or  trivial non-identifiability  problems (see Remark in Section 3.1), we will impose $p \neq 1$ and $\theta:=(\alpha,\beta)\in \Phi \subset  \R\times \R^*$, which 
implies that the Euclidean parameter $\vartheta$ will be assumed to belong to a parametric compact and convex
space 
\begin{eqnarray}\label{esp-param}
\Theta:=[\delta,1-\delta]\times\Phi\subset (0,1)\times \left\{\R\times \R^*\right\} ,
\end{eqnarray} 
where $\delta\in (0,1)$ .\\
For simplicity, we will endow the spaces $\R^s$, $s\geq 1$, with  the $\| \cdot \|_s$ norm (for clarity the dimension $s$ is  recalled in index)  defined for all $v=(v_1,\dots,v_s)$
by $\|v\|_s=\sum_{j=1}^s | v_j|$ where $| \cdot | $ denotes the absolute value.\\
We now introduce the following non-commutative  notation:  
$$
\theta\odot x:=\alpha+\beta x,\quad\quad (\theta,x)\in \Phi\times \R.
$$
Following the ideas developed by the authors mentioned above, it is  possible to 
use  the symmetry assumption made on  $f$   to identify  the
true value of the  Euclidean parameter. The idea consists in noticing  that for $\theta$
fixed in $\Phi$,   the  sample  $(Y^{\theta}_1,\dots, Y^{\theta}_n)$
 obtained by considering the so-called  {\it $\theta$-transformation}
\begin{eqnarray}\label{theta_trans}
Y^\theta_i:=Y_i-\theta \odot X_i, \quad\quad i=1,\dots,n,
\end{eqnarray} 
is distributed according to
\begin{eqnarray}\label{proj_dens}
 \Psi_{\theta}(y)=p_*\int_{\R} f(y+(\theta-\theta_*)\odot x)h(x)dx+(1-p_*) \int_{\R} f_0(y+\theta\odot x)h(x)dx,
\end{eqnarray}
where $\vartheta_*=(p_*,\alpha_*,\beta_*)\in \mathring{\Theta}$ denotes the  true value of the parameter.
Let us observe  now that  when $\theta=\theta_*$ 
\begin{eqnarray}\label{ident_formula}
\Psi_{\theta_*}(y)=p_*f(y)+(1-p_*) \int_{\R}
f_0(y+\theta_*\odot x)h(x)dx.
\end{eqnarray}
{\it Remark.} When $\theta$ is well fitted ($\theta=\theta_*$) the model associated to the $Y^\theta$ 
is very close to  the simple contamination  model (\ref{BDV06model})  studied in \cite{BDV06} or \cite{BV10}
where the location $\mu$ is known but the proportion $p$ is unknown.\\

Isolating $f$ in (\ref{ident_formula}) and replacing $\vartheta_*=(p_*,\theta_*)$ by $\vartheta=(p,\theta)$  one can define a new parametric class   of  functions  ${\cal F}_\Theta:=\left\{f_{\vartheta}:~\vartheta\in \Theta\right\}$: 
\begin{eqnarray}\label{inv_formula}
f_{\vartheta}(y)&=&\frac{1}{p}\Psi_{\theta}(y)-\frac{1-p}{p}\int_{\R} f_0(y+\theta\odot x)h(x)dx,\quad\quad(y,\vartheta)\in \R\times \Theta,
\end{eqnarray}
that satisfies  under $\vartheta=\vartheta_*$,
 \begin{eqnarray}\label{singular}
f(y)=f_{\vartheta_*}(y)=f_{\vartheta_*}(-y)=f(-y),\quad\quad y\in \R.
\end{eqnarray}
The intuition consists now in claiming that, if we make $\vartheta$ vary over $\Theta$ and that  we are able to check  that 
$f_\vartheta$ is symmetric about $0$ for a certain value of $\vartheta$ then we have reached the true value of the Euclidean parameter.
Note that in the right hand side of (\ref{inv_formula}), the second integral  term is in general unknown but can
be estimated pointwise by a standard Monte Carlo approach, see  expression  (\ref{I-estim}), or a nonparametric Monte Carlo approach, see expression (\ref{pdftilde}). The idea to check this situation, and then   to estimate $\vartheta=(p,\theta)$, is to consider a
contrast function based on the comparison between the cdf version of   $f_{\vartheta}(y)$
$$
H_1(y; \vartheta):=H_1(y;p,F_{\theta},
J_{\theta}):=\frac{1}{p}F_{\theta}(y)-\frac{1-p}{p} J_{\theta}(y),\quad\quad (y,\theta)\in \R\times\Theta,
$$
and the cdf version of $f_{\vartheta}(-y)$
$$
H_2(y; \vartheta):=H_2(y;p,F_{\theta},
J_{\theta}):=1-\frac{1}{p}F_{\theta}(-y)+\frac{1-p}{p}J_{\theta}(-y),\quad\quad (y,\theta)\in \R\times\Theta,
$$
where for all $\theta\in \Phi$,
$$J_{\theta}(y):=\int_{-\infty}^y I_{\theta}(z)dz, \quad  y\in \R,~~
\mbox{with} ~ I_{\theta}(z):=\int_{\R} f_0(z+\theta\odot x)h(x)dx,\quad z\in \R,$$
and
$$
F_{\theta}(y):=\int_{-\infty}^y f_{\theta}(z)dz,\quad\quad  y\in \R.
$$
Notice that for all $\theta$ fixed in $\Phi$, $J_\theta(\cdot)$ and $F_\theta(\cdot)$ are the cdfs associated respectively to the $\theta$-transformed known component population 
(the $Y_i$ such that $U_i=0$ in (\ref{model_arraybis})) and the $\theta$-transformed whole data.
Let us define the following  function 
\begin{eqnarray}\label{H_def}
H(y;\vartheta):=H_1(y;\vartheta)-H_2(y;\vartheta),\quad \quad (y,\vartheta)\in \R\times\Theta.
\end{eqnarray}
Notice that under $\vartheta_*$, using the symmetry of $f$,
\begin{eqnarray*}
H(y;\vartheta_*)=0,\quad \quad y\in \R.
\end{eqnarray*}
To avoid numerical integration in the approximation of an empirical contrast function based on the comparison of $H_1$ and $H_2$ over $\R$, we proceed as follows.
 Let $Q$ be an {\it instrumental } weight probability distribution  with pdf $q$ with respect to  Lebesgue measure. We suppose that $q$ is  strictly positive over $\R$ and easy to simulate. Then we  consider
\begin{eqnarray}\label{contrast_prop}
d(\vartheta):=\int_{\R} H^{2}(y,\vartheta)dQ(y),
\end{eqnarray}
where obviously $d(\vartheta)\geq 0$ for all $\vartheta\in \Theta$ and $d(\vartheta_*)=0$.
Let $(V_1,\dots,V_n)$ be an iid sample from $Q$. An empirical version  $d_n(\cdot)$ of  $d(\cdot)$ can be obtained by considering
\begin{eqnarray}\label{emp_contrast}
d_n(\vartheta):=\frac{1}{n}\sum_{i=1}^n H^2(V_i; p, \tilde
F_{n,\theta}, \hat J_{n,\theta}), \quad\quad \vartheta \in \Theta,
\end{eqnarray}
where
\begin{eqnarray}\label{J-estim}
\hat J_{n,\theta}(y):=\int_{-\infty}^y \hat I_{n,\theta}(z)dz,\quad\quad (y,\theta)\in \R\times \Phi,
\end{eqnarray}
with
\begin{eqnarray}\label{I-estim}
\hat I_{n,\theta}(z):=\frac{1}{n} \sum_{i=1}^n f_0(z+\theta\odot X_i),\quad\quad (z,\theta)\in \R\times \Phi,
\end{eqnarray}
which leads actually to the simple  expression for $\hat J_{n,\theta}(y)$
\begin{eqnarray}\label{J-estim-final}
\hat J_{n,\theta}(y):=\frac{1}{n}\sum_{i=1}^n   F_0(y+\theta\odot X_i),\quad\quad (y,\theta)\in \R\times \Phi,
\end{eqnarray}
and where  $\tilde F_{n,\theta}$ denotes a smooth version of the
empirical  cdf
$$\hat F_{n,\theta}(y):=\frac{1}{n} \sum_{i=1}^n \un_{Y_i^{\theta}\leq y},\quad\quad (y,\theta)\in \R\times \Phi,$$
defined by
\begin{eqnarray}\label{fdr_estim}
\tilde F_{n,\theta}(y):=\int_{-\infty}^{y} \hat
\Psi_{n,\theta}(t)dt,\quad\quad (y,\theta)\in \R\times \Phi,
\end{eqnarray}
where
\begin{eqnarray}\label{dens_estim}
\hat \Psi_{n,\theta}(t):=\frac{1}{n b_n}\sum_{i=1}^n
K\left(\frac{t-Y_i^{\theta}}{b_n}\right),\quad\quad (t,\theta)\in \R\times \Phi.
\end{eqnarray}
In (\ref{dens_estim}), we assume the  standard  condition insuring, for each $\theta\in \Phi$, the $L_1$ convergence of  $\hat \Psi_{n,\theta}$ towards $\Psi_{\theta}$ defined in (\ref{proj_dens})  (see Devroye \cite{D83}), namely
\begin{eqnarray}\label{Devroyecond} 
b_n\rightarrow 0, \quad \quad nb_n\rightarrow +\infty,
\end{eqnarray} 
and  $K$ is a symmetric density function. Finally we propose to estimate
$\vartheta_*$ by considering the  M-estimator
\begin{eqnarray}\label{estimateur}
\hat\vartheta_n:=(\hat p_n,\hat \theta_n)=\arg \min_{\vartheta\in
\Theta} d_n(\vartheta).
\end{eqnarray}
Once $\vartheta_*$ is estimated by $ \hat\vartheta_n$ a natural way to  estimate   $F$ and $f$ consistently is then to consider
the plug-in empirical versions of $H_1(\cdot ; \vartheta)$ and  (\ref{inv_formula}), respectively defined for all $y\in \R$ by
\begin{eqnarray}
\hat F_n(y)&:=&H_1(y;\hat p_n, \tilde F_{n,\hat \theta_n}, \tilde  J_{n,\hat\theta_n}), \label{cdf_estim}\\
\hat f_n(y)&:=&\frac{1}{\hat p_n} \hat \Psi_{n,\hat \theta_n}(y)+\frac{1-\hat p_n}{\hat p_n} \tilde I_{n,\hat \theta_n}(y),\label{pdf_estim}
\end{eqnarray}
where, for all $\theta\in \Theta$,      $\tilde I_{n,\theta}$ and $\tilde  J_{n,\theta}$ are respectively  nonparametric estimators of   $I_{\theta}$ and $J_{\theta}$  based on an iid  simulated sample $(\tilde \varepsilon_1^{[0]},\dots, \tilde \varepsilon_n^{[0]})$ from $f_0$ obtained by considering
\begin{eqnarray}
\tilde  I_{n,\theta}(t)&:=&\frac{1}{n b_n}\sum_{i=1}^n
K\left(\frac{t-(\theta\odot X_i+\tilde \varepsilon_i^{[0]})}{b_n}\right),\quad\quad (t,\theta)\in \R\times \Phi.\label{pdftilde}\\
\tilde J_{n,\theta}(y)&:=&\int_{-\infty}^{y} \tilde
I_{n,\theta}(t)dt,\quad\quad (y,\theta)\in \R\times \Phi,\label{cdftilde}
\end{eqnarray}
For convenience, the kernel used to compute (\ref{pdftilde}) will be Gaussian, {\it i.e.} $K(t)={\mathcal N}_{0,1}(t)$
where ${\mathcal N}_{m,\sigma^2}(t):=(2 \pi\sigma^2)^{-1/2} \exp((t-m)^2/2\sigma^2)$, for all $t\in \R$.
In this second plug-in step we  consider, for the  sake of  simplicity in  our proofs, the nonparametric  estimates (\ref{cdftilde}) and (\ref{pdftilde}) instead of (\ref{J-estim}) and (\ref{I-estim}). This choice  allows us
to use similar nonparametric results for both  $\hat f_{n,\theta}$ and $\tilde I_{n, \theta}$ (see the proof of Theorem \ref{Theo1} ii) and iii)), but the same results should be obtained, at the price of an aditionnal technical lemma,  by considering directly the Monte Carlo estimators (\ref{J-estim}) and (\ref{I-estim}).
\section{Identifiability and consistency}
\subsection{Identifiability}
In this section we recall  briefly why model (\ref{model}) is
identifiable under conditions similar to those established in 
\cite{BDV06} and summarized  below.
Let us define ${\mathcal{F}}_s:=\{f\in{\mathcal{F}};\int_{\mathbb{R}}|x|^sf(x)dx<+\infty\}$
for $s\geq 1$, where $\mathcal{F}$ denotes the set of even pdfs.  When  $(f,f_0)\in {\mathcal{F}}_s$ with $s\geq 2$, we denote $m:=\int_\R x^2 f(x)dx$ and $m_0:=\int_\R x^2 f_0(x)dx$.

\vspace{0.2cm}

\begin{defi}\label{defidentif}  (Identifiability). Let $(p_1,\theta_1,f_1,h_1)$
and $(p_2,\theta_2,f_2,h_2)$ denote two sets of
parameters for model (\ref{model}). The parameter in model
(\ref{model})  is said to be semiparametrically identifiable if
$$
(p_1,\theta_1,f_1(y),h_1(x))=(p_2,\theta_2,f_2(y),h_2(x)),$$
for  $\lambda^{\otimes 2}$-almost all $(x,y)\in \R^2$, whenever  we have :
\begin{eqnarray}\label{identif}
&&\left (p_{1} f_1(y-\theta_1 \odot x)+(1-p_1)
f_0(y)\right)h_1(x)\nonumber\\
&&~~~~~~~~~~~~~~~~~~~=\left (p_2f_2(y-\theta_2\odot
x)+(1-p_2) f_0(y)\right)h_2(x),
\end{eqnarray}
for  $\lambda^{\otimes 2}$-almost all $(x,y)\in \R^2$.
\end{defi}
\begin{lem}\label{Lemidentif}
If the Euclidean parameter space $\Theta$ is a subset of  $ \R\times \R^*\setminus \{0,0\}$,   $\supp(f)=\supp(f_0)=\R$,  $\supp(h)$ contains at least two intervals respectively in the neighborhood of $0$ and $+\infty$ (or $-\infty$), and the pdfs involved in model (\ref{model}) satisfy
  $(f_0,f)\in{\mathcal{F}}_3\times{\mathcal{F}}_3$, then the parameter in model (\ref{model}) is identifiable.
\end{lem}

\noindent {\it Proof.}  Integrating (\ref{identif}) with respect to $y$ over $\R$, we then obtain that $h_1(\cdot)=h_2(\cdot)$ $\lambda$-almost everywhere. Let  $h(x):=h_1(x)$ for all $x\in \supp{(h)}:=\supp{(h_1)}\cap \supp{(h_2)}$. Notice now  that, for all $x\in \supp(h)$,   (\ref{identif}) coincides with (\ref{BDV06model}) when considering
the generic location parameter $\mu$ equal to $\theta\odot x$.
In our case the first three conditional moment equations (given  $\left\{ X=x \right\}$)  associated to  (\ref{identif}) lead  to
\begin{eqnarray}\label{syst_bordes06}
\left\{
\begin{array}{ll}\displaystyle
&p_1\theta_1\odot x= p_2\theta_2\odot x,\\
&(1-p_1)m_0+p_1((\theta_1\odot x)^2+m_1) =(1-p_2)m_0+p_2((\theta_2\odot x)^2+m_2) ,\\
&p_1(3((\theta_1\odot x)m_1+(\theta_1\odot x)^3)=p_2(3((\theta_2\odot x)m_2+(\theta_2\odot x)^3).
\end{array}
\right.
\end{eqnarray}
According to \cite{BDV06}, the solutions are either,  for all $x\in \supp(h)$,
$(p_1,\theta_1\odot x)=(p_2,\theta_2\odot x)$,
which implies $(p_1,\alpha_1,\beta_1)=(p_2,\alpha_2,\beta_2)$,
or  
\begin{eqnarray}\label{syst_bordes_sol06}
\left\{
\begin{array}{lll}\displaystyle
p_2&= & {\displaystyle p_1\left(   \frac{2(\theta_1\odot x)^2}{3m_1+(\theta_1\odot x)^2-3m_0}\right)},\\
\theta_2\odot x &=&{\displaystyle \theta_1\odot x+\frac{3m_1-(\theta_1\odot x)^2-3m_0}{2 \theta_1\odot x}},\\
m_2&=&{\displaystyle m_1+\frac{(m_1+ (\theta_1\odot x)^2-m_0)(3m_1+ (\theta_1\odot x)^2-3m_0 )}{4(\theta_1\odot x)^2}}.
\end{array}
\right.
\end{eqnarray}
Suppose that $\beta_1\neq 0$ and take  the limit as $x\rightarrow +\infty$ in the first row of (\ref{syst_bordes_sol06}). We then necessarily
obtain that $p_2=2p_1$ which is only compatible, when we take the limit as $x\rightarrow 0$,  with $m_1=m_0$.
Hence if  $m_1\neq m_0$ model (\ref{model}) is always   identifiable. If we suppose $m_1=m_0$, the second row of (\ref{syst_bordes_sol06})
leads to $\theta_2\odot x=(\theta_1\odot x)/2$. If we introduce  this last relation in the third row of $(\ref{syst_bordes_sol06})$ we obtain
$$
m_2-m_1=\frac{1}{4} (\theta_1\odot x)^2, \quad\quad x\in \R,
$$
which is impossible when $x\rightarrow +\infty$ and thus provides us the global identifiability of model (\ref{model}).
\hfill $\square$\\

\noindent \textit{Remark. }  In Lemma \ref{Lemidentif} we have considered for simplicity the case where the slope parameter  $\beta$ is   supposed to be different away from zero. Actually  this condition can be technically   relaxed if we allow $\theta$ to be equal to  $(\alpha,0)$ with  $\alpha\neq 0$.  In fact, considering 
the  first row of (\ref{syst_bordes_sol06}) and taking the limit  as $x\rightarrow +\infty$, we obtain $\beta_2=\beta_1=0$. To conclude, it is then enough to integrate (\ref{identif}) with respect to $x$ over $\R$ which leads to discuss the same condition as  in \cite{BDV06},  p. 735 expression (3). Then Proposition 2 in  \cite{BDV06}   provides  an {\it almost everywhere}-type identifiability result  which unfortunately cannot be strictly  compared  to  the result stated in  Lemma 3.1. For this  reason we decided to reject  $\theta=(\alpha,0)$,  $\alpha\in \R^*$, from the sub-parametric space $\Phi$ , see (\ref{esp-param}).

 \subsection{Assumptions and statistical complexity} In the following we provide some general conditions that allow us to control the statistical complexity of our model and that insure the validity of basic asymptotic results.\\
 
\noindent {\bf Regularity conditions (R).}
 \begin{enumerate}[i)]
 
 \item The pdfs $f$ and $f_0$ are   strictly positive over $\R$ and  belong to ${\mathcal F}_3$.
 
 \item The  pdfs $f$ and $f_0$ are  twice differentiable over $\R$ with $\|f^{(j)}\|_\infty<\infty$ and $\|f_0^{(j)}\|_\infty<\infty$,
 where  $f^{(j)}$ and $f_0^{(j)}$ denote respectively  the $j$-th order derivatives  of $f$ and $f_0$, for $j=1,2$.
 
 \item The pdf $h$ satisfies $\int_\R |x|^2h(x)dx<\infty$.
 
\item For $i=0$ or $i=2$,
 $$\int_{\R^2} |x|^{i}|F_0(y+\theta_*\odot x)-F_0(y-\theta_*\odot x)|h(x)dxdy<\infty, $$
 and for $i=1$ or $i=3$,  and all $u\in \R$, 
 $$
 \lim_{ y\rightarrow \pm \infty} y^{i} (F_0(y+u)-F_0(y-u))=0.
 $$
 
 \item There exist  two collections of functions $\left\{\ell_{i,j}\right\}_{0\leq i\leq j\leq 2}$ and $\left\{\ell^0_{i,j}\right\}_{0\leq i\leq j\leq 2}$ belonging  to $L_1(\R^2)$ and   such that,  for all $(x,y)\in \R^2$ and all $\theta\in \Theta$
 $$|x^{i}   f^{(j)}(y+(\theta-\theta^*)\odot x)| h(x)\leq  \ell_{i,j}(x,y),$$
 and
 $$|x^{i}   f_0^{(j)}(y+\theta\odot x)| h(x)\leq  \ell^0_{i,j}(x,y).$$

 \end{enumerate}

\noindent  For all  $z\in \C$, let  $\breve{z}$ and $\Im(z)$  the conjugate and   imaginary part  of $z$, respectively. We will also denote     $\bar f$,  $\bar f_0$ the Fourier transforms of $f$, $f_0$,  and define for all $\kappa=(\kappa_1,\kappa_2)\in \R^2$, $\nu_{\kappa}(t):=e^{i t\kappa_1} \bar h(\kappa_2t)$,   where $\bar h$ denotes the Fourier transform of $h$.\\

The following conditions mainly insure the contrast property for the function $d$ defined in (\ref{contrast_prop}). We point out  that these conditions are not equivalent, as is the case in \cite{BV10} p. 25, to those established to prove the identifiability property in Lemma \ref{Lemidentif}. Loosely speaking the reason of this difference is due to the $\theta$-transformation that reduces the Euclidean parameter estimation problem to the analysis of a collection of one-dimensional data, {\it i.e.} the $Y^{\theta}_i$ with $\theta\in \Phi$, when the proof of Lemma \ref{Lemidentif} uses strongly the bivariate  structure of the original data.\\

 \noindent{\bf Contrast condition (C).} 
 \begin{enumerate}[i)]
\item  The three  first moments of $X$ satisfy $4E(X)^3+3E(X)E(X^2)+E(X^3)\neq 0.$
\item The set of parameters $\vartheta=(p,\theta)=(p,\alpha,\beta)$ with $p\neq p_*$ that satisfies
 \begin{eqnarray}\label{Fourier_eq}
p_*\Im(\nu_{\theta-\theta_*}(t))\bar f(t)= (p_*-p)\Im(\nu_{\theta}(t))\bar f_0(t),\quad\quad t\in \R,
\end{eqnarray}
 is empty or does not belong to the parametric space $\Theta$.

\item The second order moments of $f$ and $f_0$, respectively denoted $m$ and $m_0$, are supposed to satisfy
$$
m\neq m_0+\frac{\alpha_*^3+3\alpha_*^2\beta_* E(X)+3\alpha_* \beta_*^2 E(X^2)+\beta_*^3 E(X^3)}{3(\alpha_*+\beta_* E(X))}.
$$ 
\end{enumerate}
{\it Remark.} Point out that condition C ii), which is necessary  to prove that $d$ is a contrast  function over $\Theta$, cannot be simplified
without more information on $f$, $f_0$ and $h$.  We suggest, in  the spirit of conditions C1 and C2  in \cite{HH11},  to  
consider the  sufficient and more intuitive   \textit {regularity comparison}-type criterion for C ii)
\begin{eqnarray}\label{critere}
\forall \theta\in \Phi, \quad \left |\frac{\Im(\nu_{\theta-\theta_*}(t))} {\Im(\nu_{\theta}(t))} \frac{\bar f(t)}{\bar f_0(t)}\right |\longrightarrow +\infty ~~ \mbox{or} ~~0,~~ \mbox{as}~~ t\rightarrow +\infty,
\end{eqnarray}
which is valid since,  according to (\ref{Fourier_eq}), the term on left  hand side of  (\ref{critere}) is equal  to $|p-p_*|/p_*\in[|p-p_*|, 1/\delta]$  which is in contradiction with (\ref{critere}).
 However  condition (\ref{Fourier_eq}) can directly  be  discussed in the  Gaussian case
as  done in the appendix, Section \ref{ident_gauss_sect}.  We prove in particular that there exist sometimes spurious solutions satisfying
(\ref{Fourier_eq}) that have to be removed from the parametric space so they are not detected by our estimation algorithm as shown in Fig.  \ref{fig:fig_shap12}.
 \\

 \noindent \textbf{Kernel and Bandwidth conditions (K). }
\begin{enumerate}[i) ]
\item The even kernel density function $K$ is bounded, uniformly continuous, square integrable, of bounded variation and has second order moment.
\item The bandwidth $b_n$ satisfies $b_n\searrow 0$, $nb_n\to+\infty$ and $\sqrt{n}b_n^2=o(1)$.
\end{enumerate}

\begin{lem} \label{LemTech}
\begin{enumerate}[(i)]
\item Under conditions (R)  the function  $d$ is Lipschitz  over  $\Theta$.
\item  Under conditions (C) i) and ii) the function  $d$ is a contrast function, {\it i.e.} for all $\vartheta \in \Theta$, $d(\vartheta)\geq 0$ and $d(\vartheta)=0$ if and only if $\vartheta=\vartheta_*$.
\item Under condition (C) iii) we have 
$$
\ddot
d(\vartheta_*)=2\int_\mathbb{R}\Dot{H}(y,\vartheta_*)\Dot{H}^T(y,\vartheta_*)dQ(y)>0.
$$
\item Under conditions (R) and (K),  for any $\gamma>0$, $d_n$ converges  to $d$ almost surely   with the rate    $$\sup_{\vartheta\in \Theta}|d_n(\vartheta)-d(\vartheta)|=o_{a.s.}(n^{-1/2+\gamma}).$$
\end{enumerate}
\end{lem}

\noindent{\it Remark .}  There exists a simple consistent  method  to select, in the $L_1(\R^2)$ sense (recall that our nonparametric consistency results are established for this norm),  the best   estimator in case of multiple minima of $d_n$ (which should  make suspect that condition (C) is violated).
Suppose that,  for $n$ fixed in $\N^*$,  there exists a finite collection of local minima of $d_n$, denoted by $\hat \vartheta_n^{[i]}=\left (\hat p^{[i]}_n, \theta_n^{[i]}\right)$  with  $1 \leq i\leq S<\infty$. Then we propose to retain a  $\hat \vartheta_n$ (in practice unique) satisfying  
$$
\hat \vartheta_n=\hat \vartheta_n^{[i_*]}, ~~ \mbox{where}~~i_*=  \arg \min_{1\leq i \leq S} \left\|\hat g_n-\hat g_{\hat \vartheta_n^{[i]}}\right\|_{L_1},
$$
and where for all $1\leq i\leq S$  ,  $\hat g_{\hat \vartheta_n^{[i]}}$ is the plug-in posterior estimator of $g$ defined by 
\begin{eqnarray}
\hat g_{\hat \vartheta_n^{[i]}}=\hat p_n^{[i]} \hat f_{\vartheta_n^{[i]}}+(1-\hat p_n^{[i]}) f_0,
\end{eqnarray}
where $\hat f_{\vartheta_n^{[i]}}$ corresponds to $f_n$ defined in (\ref{pdf_estim}), when $\hat \vartheta_n=\hat\vartheta_n^{[i]}$. Proceeding in that way,  we clearly support  the Euclidean parameter estimate that better fit the dataset, this approach being asymptotically consistent as long as the  the model is identifiable.\\

\noindent\textsc{Proof. } i) From boundedness and  the uniform  Lipschitz property of  $H(\cdot, \vartheta)$, along with the  integrability  and the integrable Lipschitz property of $f_\theta(\cdot)$  proved in Sections \ref{bound}, \ref{integ_lipf} and \ref{integ_lipH}, there exists  a nonnegative constant   $c$ such that  for all $(\vartheta,\vartheta ')\in \Theta^2$
\begin{eqnarray*}
&&\left|\int_\R H^2(y,\vartheta) dQ(y)-\int_\R H^2(y,\vartheta') dQ(y)\right |\\
&\leq&\int_\R \left | H(y,\vartheta)+H(y,\vartheta') \right |  \left | H(y,\vartheta)-H(y,\vartheta') \right |q(y)dy\\
&\leq&c\|\vartheta-\vartheta'\|_3,
\end{eqnarray*}
which concludes the proof of i).

ii) To clarify the similarity between the semiparametric contamination
model  (\ref{BDV06model}) studied in \cite{BDV06} and the contaminated
regression model (\ref{model}), we can  say that
$f_\theta(\cdot)$ plays the role of $g(\cdot-\mu)$ and that
$I_\theta(\cdot)$ plays the role of $f_0(\cdot-\mu)$.\\

If $\vartheta=\vartheta_*$ then $d(\vartheta)=0$. To prove the
converse we notice that $d(\vartheta)=0$ implies, since $H_1(\cdot,\vartheta)$ and $H_2(\cdot,\vartheta)$
are continuous and $q>0$ over $\R$, that 
$H_1(\cdot;\vartheta)=H_2(\cdot;\vartheta)$ which leads,   for almost all $y\in \R$, to 
\begin{eqnarray}\label{posit_contrast}
f_\theta(y)-(1-p)I_\theta(y)=f_\theta(-y)-(1-p)I_\theta(-y).
\end{eqnarray}
Using formula (\ref{proj_dens}), we obtain
\begin{eqnarray*}
&& p_*\int_\R f(y+(\theta-\theta_*)\odot x)h(x)dx+(p-p_*)I_\theta(y)\\
&~&~~~~~~~~~~=p_0\int_\R
f(-y+(\theta-\theta_*)x)h(x)dx+(p-p_*)I_\theta(-y),\quad \quad y\in
\R.
\end{eqnarray*}
Considering the Fourier transform of the previous equality,  using
 Fubini's Theorem, and noticing that $\bar f$
and $\bar f_0$ are real-valued functions,  it follows  that
\begin{eqnarray*}
&&p_*   e^{-it(\alpha-\alpha_*)}\bar f(t)\breve {\bar h}((\beta-\beta_*)t)+(p-p_*) e^{-it \alpha} \bar f_0(t) \breve{\bar
h}(\beta t) \\
&&~~~~~~~~~~~~~~=p_*   e^{it(\alpha-\alpha_*)}\bar f(t)\bar h((\beta-\beta_*)t)+(p-p_*) e^{it \alpha} \bar f_0(t) \bar
h(\beta t),\quad\quad t\in \R.
\end{eqnarray*}
Using the notation introduced for the writing of condition (C),
the previous equation becomes (\ref{Fourier_eq}).\\

Suppose that $p=p_*$  and take the first and third order derivative of (\ref{Fourier_eq}) at point $t=0$. We then obtain $\alpha-\alpha_*+(\beta-\beta_*)E(X)=0$ and 
\begin{eqnarray*}
&&3m[\alpha-\alpha_*+(\beta-\beta_*)E(X)]+(\alpha-\alpha_*)^3+3(\alpha-\alpha_*)^2(\beta-\beta_*)E(X)\\
&&~~~~~~~~~~~+3(\alpha-\alpha_*)(\beta-\beta_*)^2E(X^2)+(\beta-\beta_*)^3E(X^3)=0,
\end{eqnarray*}
which naturally leads to 
\begin{eqnarray}
(\beta-\beta_*)(4E(X)^3+3E(X)E(X^2)+E(X^3))=0,
\end{eqnarray}
and thus implies that $\theta=\theta_*$ if  $4E(X)^3+3E(X)E(X^2)+E(X^3)\neq0$.\\

Suppose now that $p\neq p^*$, then  condition (C) ii)  requires that  $\vartheta=\vartheta_*$.\\

iii) First we have
\begin{eqnarray}\label{hess_theta*}
\ddot d(\vartheta_*)&=&2\int_\R\left( \ddot H(y,\vartheta_*)H(y,\vartheta_*)+ \dot H(y,\vartheta_*)\dot H^T(y,\vartheta_*)\right ) q(y)dy\nonumber\\
&=&2\int_\R \dot H(y,\vartheta_*)\dot H^T(y,\vartheta_*) q(y) dy,
\end{eqnarray}
according to (\ref{ident_formula}) and the fact that  $H(\cdot,\vartheta_*)=0$ on $\R$. Let $ v$ be a vector  in $\R^3$. We have 
\begin{eqnarray}
v^T \ddot d(\vartheta_*)v= 2\int_\R \left( v^T  \dot H(y,\vartheta_*)\right)^2q(y) dy\geq 0.
\end{eqnarray}
It follows that $\ddot d(\vartheta_*)$ is a positive $3\times 3$ real valued matrix. Let us show that it is also definite. If $v\in \R^3$ is a non-null column vector such that 
$v^T \ddot d(\theta_*) v=0$, then $v^T\dot H(y,\vartheta_*)=0$  for almost all  $y\in\R$.  According to (\ref{deriv H}) in the appendix, we have to discuss
the proportionality of $f$ and $F_0^{\theta_*}(\cdot)+F_0^{\theta_*}(-\cdot)-1$. Because $f_0$ is an even density, we have from   Fubini's theorem
\begin{eqnarray}\label{f_isol}
f(y)=\frac{\int_\R [F_0(y+\theta_*\odot x)-F_0(y-\theta_*\odot x)]h(x)dx}{\int_{\R^2} [F_0(y+\theta_*\odot x)-F_0(y-\theta_*\odot x)]h(x)dxdy}.
\end{eqnarray}
Using  integration by parts and assumption (R) iv), the   denominator of  the right hand side  of  (\ref{f_isol}) can be expressed as follows
\begin{eqnarray*}
&&\int_{\R^2} [F_0(y+\theta_*\odot x)-F_0(y-\theta_*\odot x)]h(x)dxdy\\
&&~~~~~~~~=\int_\R \left\{\left[  y (F_0(y+\theta_*\odot x)-F_0(y-\theta_*\odot x)   \right]_{-\infty}^\infty  \right\}h(x)dx\\
&&~~~~~~~~~-\int_\R \int_\R y(f_0(y+\theta_*\odot x)-f_0(y-\theta_*\odot x) )dy h(x)dx\\
&&~~~~~~~~=2\int_\R (\alpha_*+\beta_* x)h(x)dx=2(\alpha_*+\beta_* E(X)).
\end{eqnarray*}
If we calculate now the second-order moment of $f$ we obtain
\begin{eqnarray*}
m&:=&\int_\R y^2f(y)dy=\frac{\int_\R x^2[F_0(y+\theta_*\odot x)-F_0(y-\theta_*\odot x) ]h(x)dx}{2(\alpha_*+\beta_* E(X))}.
\end{eqnarray*}
Using integration by parts and assumption (R) iv), the numerator of the  right hand-side of (\ref{f_isol})   can be expressed as follows
\begin{eqnarray*}
&&\int_{\R^2}y^2 [F_0(y+\theta_*\odot x)-F_0(y-\theta_*\odot x)]h(x)dxdy\\
&&~~~~~~~~=\int_\R \left\{\left[  \frac{y^3}{3} (F_0(y+\theta_*\odot x)-F_0(y-\theta_*\odot x)   \right]_{-\infty}^\infty  \right\}h(x)dx\\
&&~~~~~~~~~-\int_{\R^2}\frac{y^3}{3}(f_0(y+\theta_*\odot x)-f_0(y-\theta_*\odot x) )dy h(x)dx\\
&&~~~~~~~~=2\int_{\R^2} \frac{3u^2\theta\odot x+(\theta\odot x)^3}{3}f_0(u)du h(x)dx\\
&&~~~~~~~~=2m_0(\alpha_*+\beta_* E(X))+\frac{2}{3}(\alpha_*^3+3\alpha_*^2\beta_* E(X)+3\alpha_* \beta_*^2 E(X^2)+\beta_*^3 E(X^3)),
\end{eqnarray*}
which leads to a contradiction if (C) iii) is assumed.\\

iv) This proof, which is a tricky   generalization  of the proof of   Lemma 3.2 iii) given  in \cite{BV10}, is relegated to the appendix for convenience, see Section \ref{unif_rated}.\hfill $\square$

\begin{thm} \label{Theo1}
\begin{enumerate}[i)]
\item If assumptions (R), (C) and (K) are satisfied then 
$$
\|\hat \vartheta_n-\vartheta_*\|_3=o_{a.s.}(n^{-1/4+\gamma}),\quad \quad \gamma>0.
$$
\item The estimator $\hat f_n$ of $f$ defined in (\ref{pdf_estim}) converges almost surely in the $L_1$ sense 
if  $n^{-1/4+\gamma}/b_n\rightarrow 0$, for all $\gamma>0$.
\item For any $\gamma>0$, the estimator $\hat F_n$ of $F$ defined in (\ref{cdf_estim}) converges uniformly at the following  almost sure rate
\begin{eqnarray}\label{rateF}
\|\hat F_n-F \|_\infty=O_{a.s.}(n^{-1/4+\gamma}/b_n)+O_{a.s.}(b_n^2), \quad\quad  \gamma>0.
\end{eqnarray}
The above rate is optimized by considering  $b_n=n^{-1/12}$, which choice provides the rate of convergence
 $O_{a.s.}(n^{-1/6+\gamma})$, for all  $\gamma>0$.
\end{enumerate}
\end{thm}
{\it Comment.} Points ii) and iii) reveal  the intuitive idea  that the bandwidth $b_n$ must not decrease too fast in order
to allow the appropriate  positionning of the plug-in-centered data in the expression of $\hat \Psi_{n, \hat \theta_n}$. In fact the  $Y^{\hat\theta_n}_i$ need to be sufficiently close  to the $Y^{\theta_*}_i$, and $b_n$ not too small (the smaller $b_n$ is  the more we ``freeze" the kernel estimator),
if we want a good agreement   between $\hat \Psi_{n, \hat \theta_n}$ and  $\hat \Psi_{n, \theta_*}$ which are  known to converge to the true $\Psi_{\theta_*}$ involved in  expression (\ref{ident_formula}).\\

\noindent\textsc{Proof.} 
i) The proof follows entirely  the proof of  Theorem 3.1 in \cite{BV10} and uses the technical results proved in 
Lemma \ref{LemTech}.

ii) Consider the following decomposition:
\begin{eqnarray}\label{decompL1}
|\hat f_n-f|
&=&\left |\left[\frac{1}{\hat p_n} \hat \Psi_{n,\hat \theta_n}-\frac{1}{p_*} \Psi_{\theta_*}\right]+\left[\frac{1-\hat p_n}{\hat p_n} \tilde I_{n,\hat \theta_n}-\frac{1- p_*}{ p_*}  I_{\theta_*} \right]\right|\nonumber\\
&=&\left|\left[\frac{1}{\hat p_n} (\hat \Psi_{n,\hat \theta_n}- \hat \Psi_{n,\theta_*})+\frac{1}{\hat p_n} \hat \Psi_{n,\theta_*}-\frac{1}{p_*} \Psi_{\theta_*}\right]\right. \nonumber\\
&&\left.+\left[\frac{1-\hat p_n}{\hat p_n} (\tilde I_{n,\hat \theta_n}- \tilde I_{n,\theta_*})+\frac{1- \hat p_n}{ \hat p_n}  \tilde I_{n,\theta_*}-I_{\theta_*} \right]\right|\nonumber\\
&\leq& \frac{1}{\hat p_n}\left (|\hat \Psi_{n,\hat \theta_n}- \hat \Psi_{n,\theta_*}|+|\hat \Psi_{n,\theta_*}- \Psi_{\theta_*}|\right)+\Psi_{\theta_*} \left|\frac{1}{\hat p_n}-\frac{1}{p_*}\right|\nonumber\\
&&+ \frac{1-\hat p_n}{\hat p_n} \left (|\tilde I_{n,\hat \theta_n}- \tilde I_{n,\theta_*}|+|\tilde I_{n,\theta_*}- I_{\theta_*}|\right)\nonumber\\
&&+I_{\theta_*} \left|\frac{1-\hat p_n}{\hat p_n}-\frac{1-p_*}{p_*}\right|.
\end{eqnarray}
It is now enough to study the behavior of $|\hat \Psi_{n,\hat \theta_n}- \hat \Psi_{n,\theta_*}|$ and  $|\tilde I_{n,\hat \theta_n}- \tilde I_{n,\theta_*}|$. For all $t\in \R$, we have 
\begin{eqnarray}\label{approxf}
|\hat \Psi_{n,\hat \theta_n}(t)- \hat \Psi_{n,\theta_*}(t)|&\leq& \frac{1}{nb_n}\sum_{k=1}^n \left | K\left (\frac{t-Y_i^{\hat \theta_n}}{b_n}\right) - K\left (\frac{t-Y_i^{\theta_*}}{b_n}\right) \right|.
\end{eqnarray} 
Consider $K$ a centered normalized  gaussian kernel. We propose to study in a generic way the difference of kernels involved in the right hand side of the above expression.  For all $(w, z)\in \R^2$, and letting  $h:=(z-w)/b$,  we write the second-order Taylor expansion with integral remaining term:
\begin{eqnarray*}
K\left (\frac{t-w}{b}\right)-K\left (\frac{t-z}{b}\right)=h \dot K\left (\frac{t-z}{b}\right)+\frac{h^2}{2}\int_0^1 (1-u) \ddot K\left (\frac{t-m_u}{b}\right)du
\end{eqnarray*}
where $m_u:=(1-u)z+uw$.
Noticing that 
\begin{eqnarray*}
\dot K\left (\frac{t-z}{b}\right)=(t-z){\mathcal N}_{z,b^2}(t), \quad\mbox{and}\quad \ddot K\left (\frac{t-m_u}{b}\right)=\left(-b+\frac{(t-m_u)^2}{b}\right){\mathcal N}_{m_u,b^2}(t)dt,
\end{eqnarray*}
it thus follows  that 
\begin{eqnarray}\label{approxker}
&&\frac{1}{b}\int_{\R}\left|K\left (\frac{t-w}{b}\right)-K\left (\frac{t-z}{b}\right)  \right|dt\nonumber\\
&&\leq h\int_\R |t-z| \frac{{\mathcal N}_{z,b^2}(t)}{b}dt +\frac{h^2}{2}\int_\R (1+\frac{(t-m_u)^2}{b^2}){\mathcal N}_{m_u,b^2}(t)dt\nonumber\\
&&\leq -h \sqrt{\frac{2}{\pi }}\int_{0}^\infty - \frac{t}{b^2} \exp(-\frac{t^2}{2b^2})dt + h^2\nonumber\\
&& \leq \sqrt{\frac{2}{\pi }}h+h^2.
\end{eqnarray}
Replacing   $w$, $z$ respectively by the $Y_i^{\theta_n}$ and $Y_i^{\theta_*}$,  and $b$ by $b_n$ in (\ref{approxker}) we then obtain from (\ref{approxf}) 
the following bound for the $L_1$ error:
\begin{eqnarray}\label{approxkey}
\|\hat \Psi_{n,\hat \theta_n}(t)- \hat \Psi_{n,\theta_*}(t)\|_{L_1}\leq \frac{C\|\theta_n-\theta_*\|_2}{b_n}\times \frac{1}{n}\sum_{k=1}^n 
\left(|X_i|+|X_i|^2\right), 
\end{eqnarray}
the same kind of bound being available  for $\|\tilde  I_{n,\hat \theta_n}(t)- \tilde  I_{n,\theta_*}(t)\|_{L_1}$. In conclusion, according  to the decomposition (\ref{decompL1}), point i) of Theorem \ref{Theo1}, the respective $L_1$ $a.s.$ convergence of $\hat \Psi_{n,\theta_*}$ and $\tilde I_{n,\theta_*}$ towards $\Psi_{\theta_*}$ and $I_{\theta_*}$ under (\ref{Devroyecond}), we get from (\ref{approxkey}) and the strong law of large numbers that $\|\hat f_n-f\|_{L_1}\rightarrow 0$ almost surely as $n\rightarrow \infty$ whenever $n^{-1/4+\gamma}/ b_n=o(1)$.\\

iii) The proof  uses an integrated version of decomposition (\ref{decompL1}) and the fact that, for all $y\in \R$,  the approximation $|\hat F_{n,\hat \theta_n}- \hat F_{n,\theta_*}|(y)$  is  controlled by 
\begin{eqnarray}\label{approxkeybis}
|\hat F_{n,\hat \theta_n}- \hat F_{n,\theta_*}|(y)&=&\left|\int_{-\infty}^y \frac{1}{n}\sum_{i=1}^n K\left (\frac{t-Y_i^{\hat \theta_n}}{b_n}\right) - K\left (\frac{t-Y_i^{\theta_*}}{b_n}\right)  \right|dt \nonumber\\
&\leq&\frac{1}{n}\sum_{i=1}^n\int_\R \left |K\left (\frac{t-Y_i^{\hat \theta_n}}{b_n}\right) - K\left (\frac{t-Y_i^{\theta_*}}{b_n}\right) \right|dt\nonumber\\
&\leq&\frac{C\|\theta_n-\theta_*\|_2}{b_n}\times \frac{1}{n}\sum_{k=1}^n 
\left(|X_i|+|X_i|^2\right),
\end{eqnarray}
the last term in the right hand side of above inequality being independent from $y$.  The same bound holds for  $|\tilde I_{n,\hat \theta_n}- \tilde I_{n,\theta_*}|(y)$ by an identical argument. To conclude, it is enough to use (\ref{approxkeybis}) and  Corollary 1 p. 766 in \cite{SW86} which allows us to control the terms  $\|\hat F_{n,\theta_*} -F_{\theta_*}\|_\infty$ and $\|\tilde I_{n,\theta_*} -I_{\theta_*}\|_\infty$, to obtain (\ref{rateF}). The rate on  right hand side of (\ref{rateF}) is optimized by considering 
$b_n=n^{-1/12}$ which  then turns into $O_{a.s.}(n^{-1/6+\gamma})$, for all $\gamma>0$.

\hfill $\square$

\section{Numerical experiments}
\subsection{Role of the $\theta$-transformation}
\begin{figure}[h!]
\includegraphics[width=3.9cm]{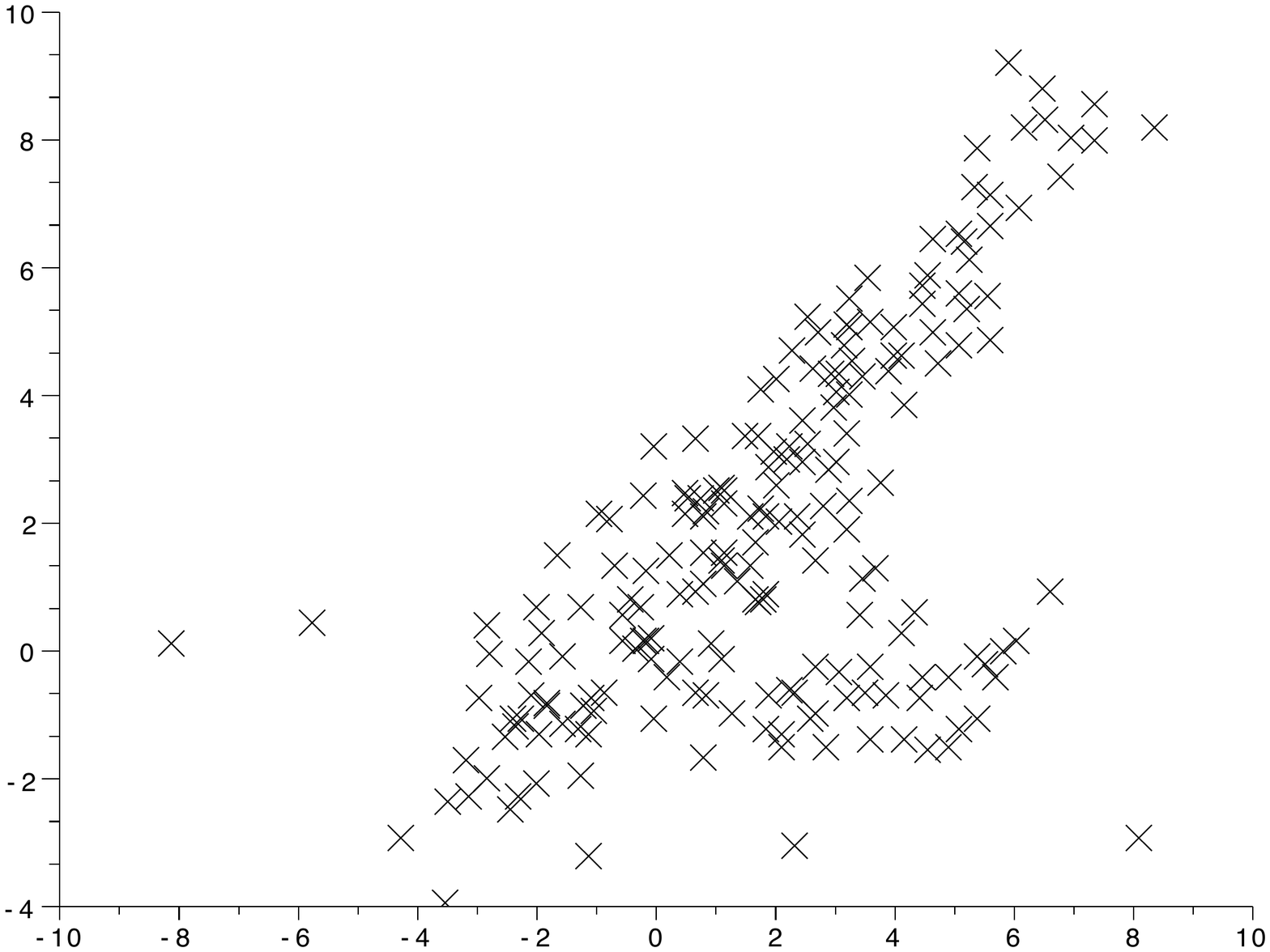}
\includegraphics[width=3.9cm]{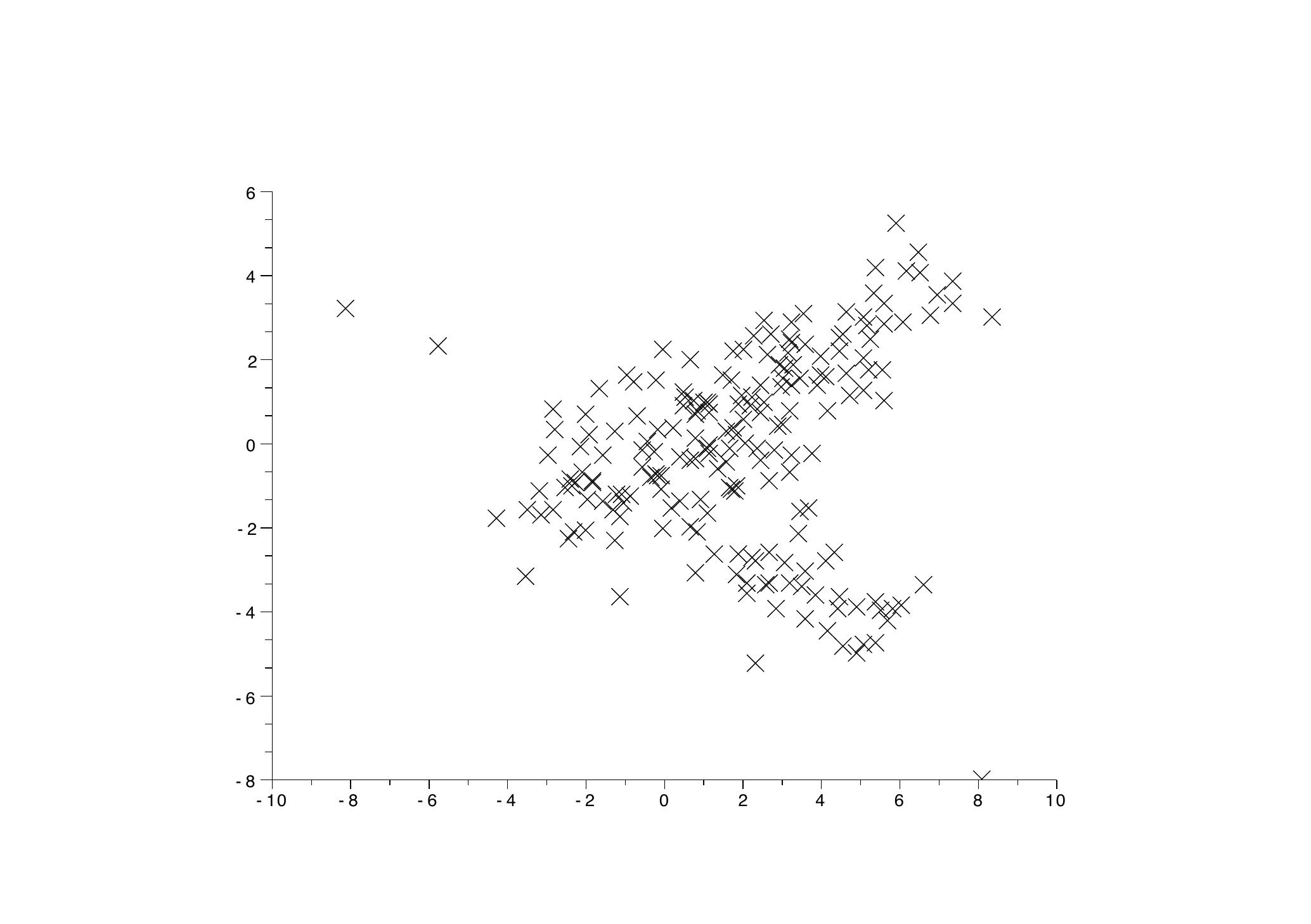}
\includegraphics[width=3.9cm]{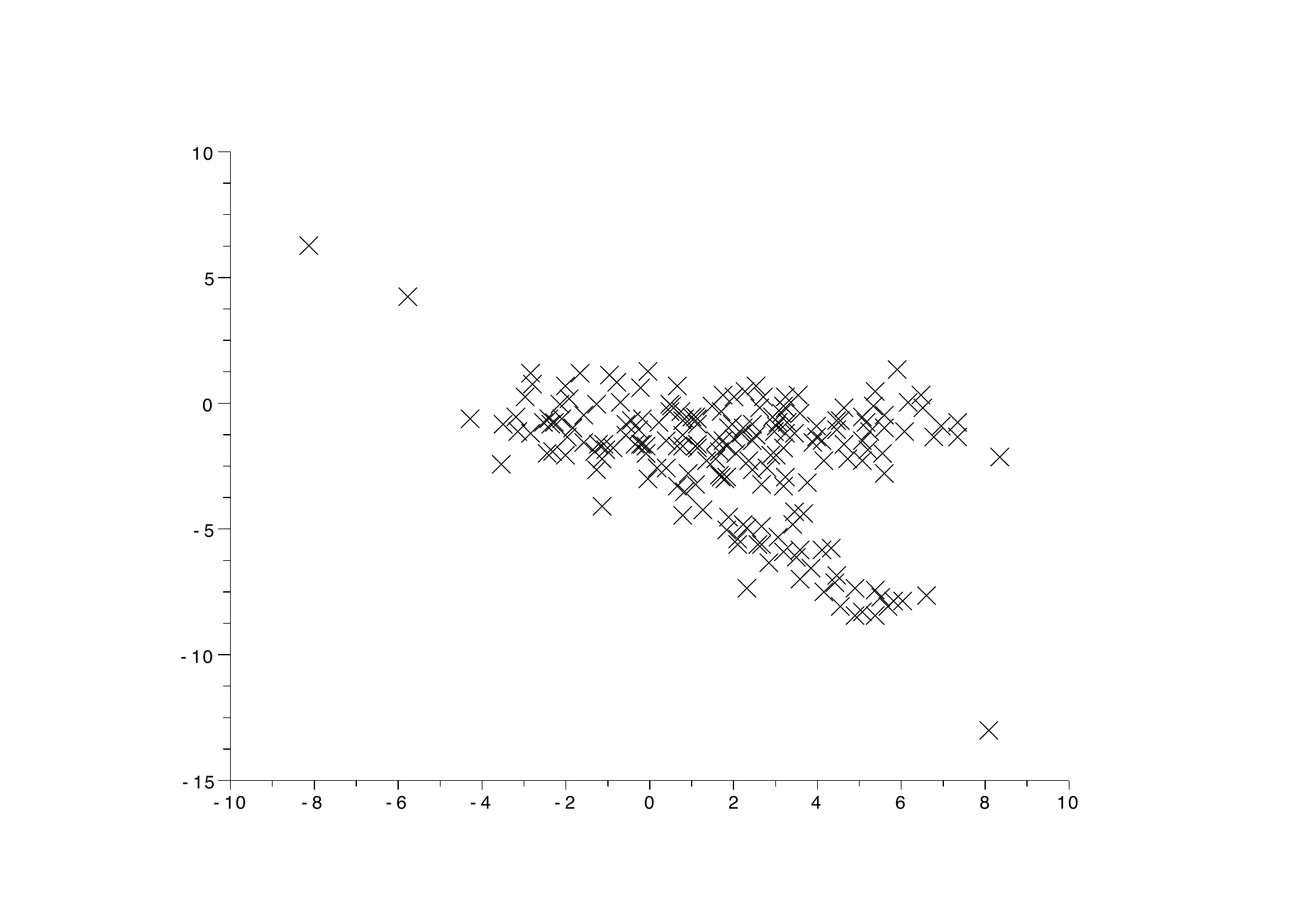}\\
\includegraphics[width=3.9cm]{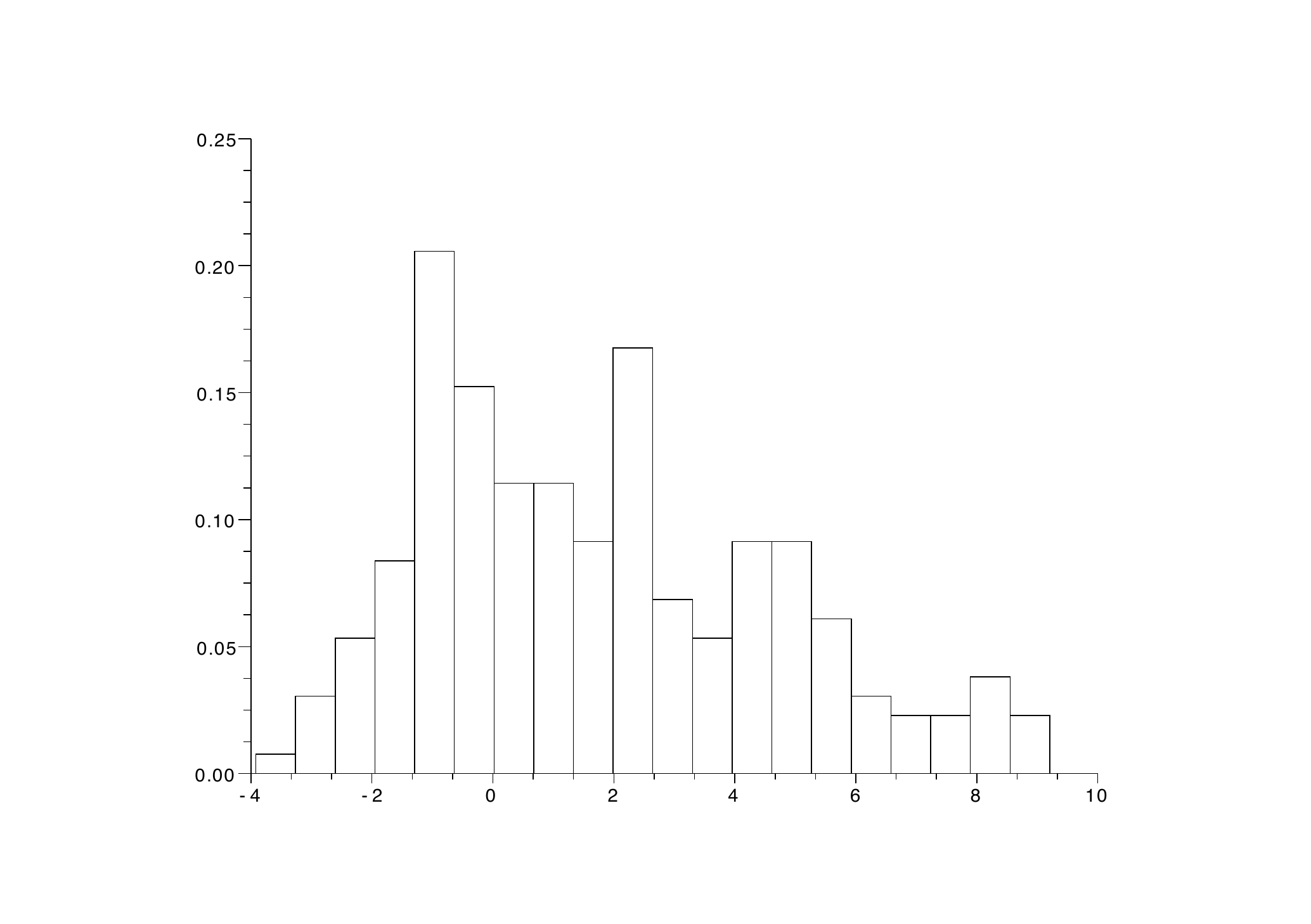}
\includegraphics[width=3.9cm]{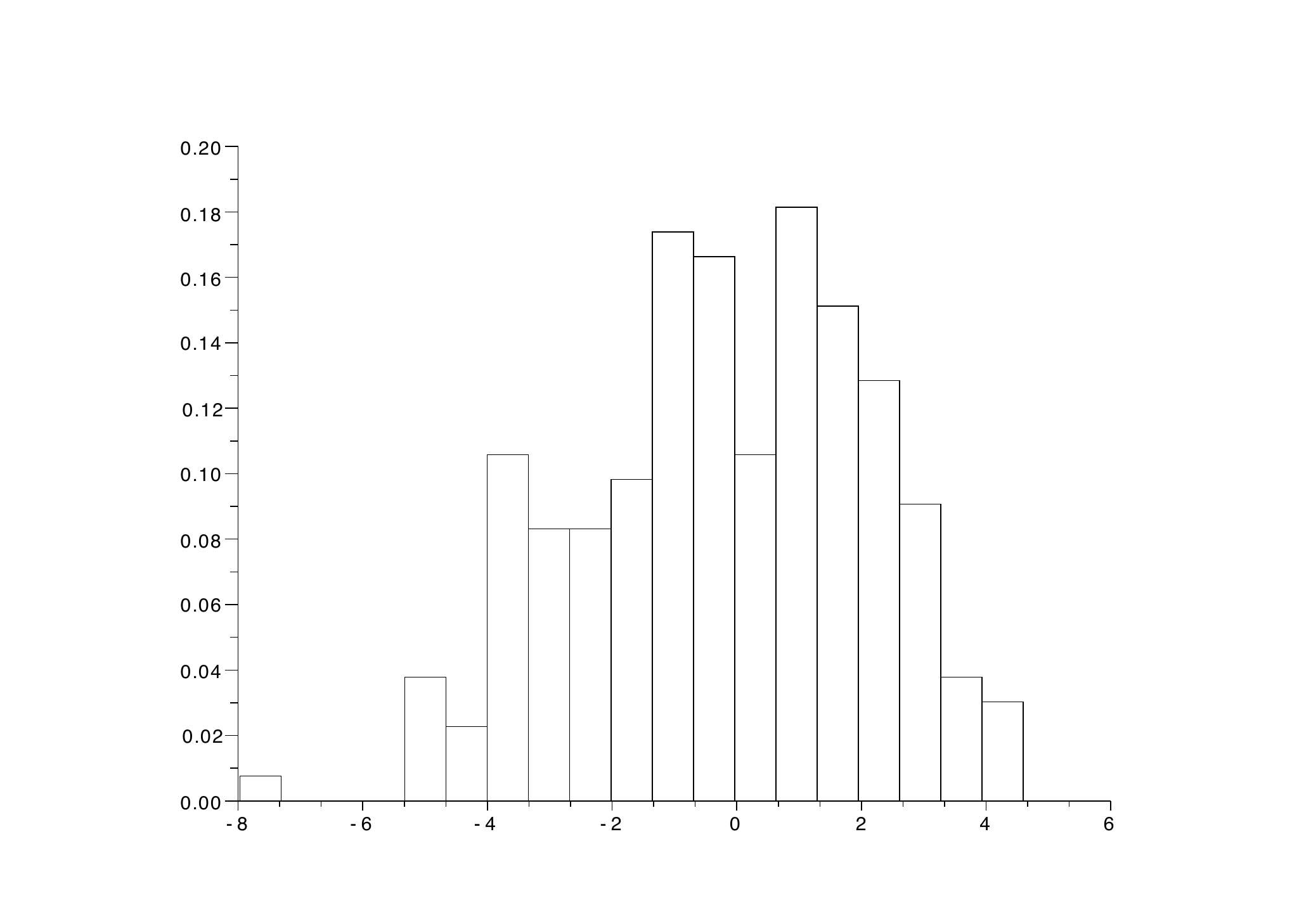}
\includegraphics[width=3.9cm]{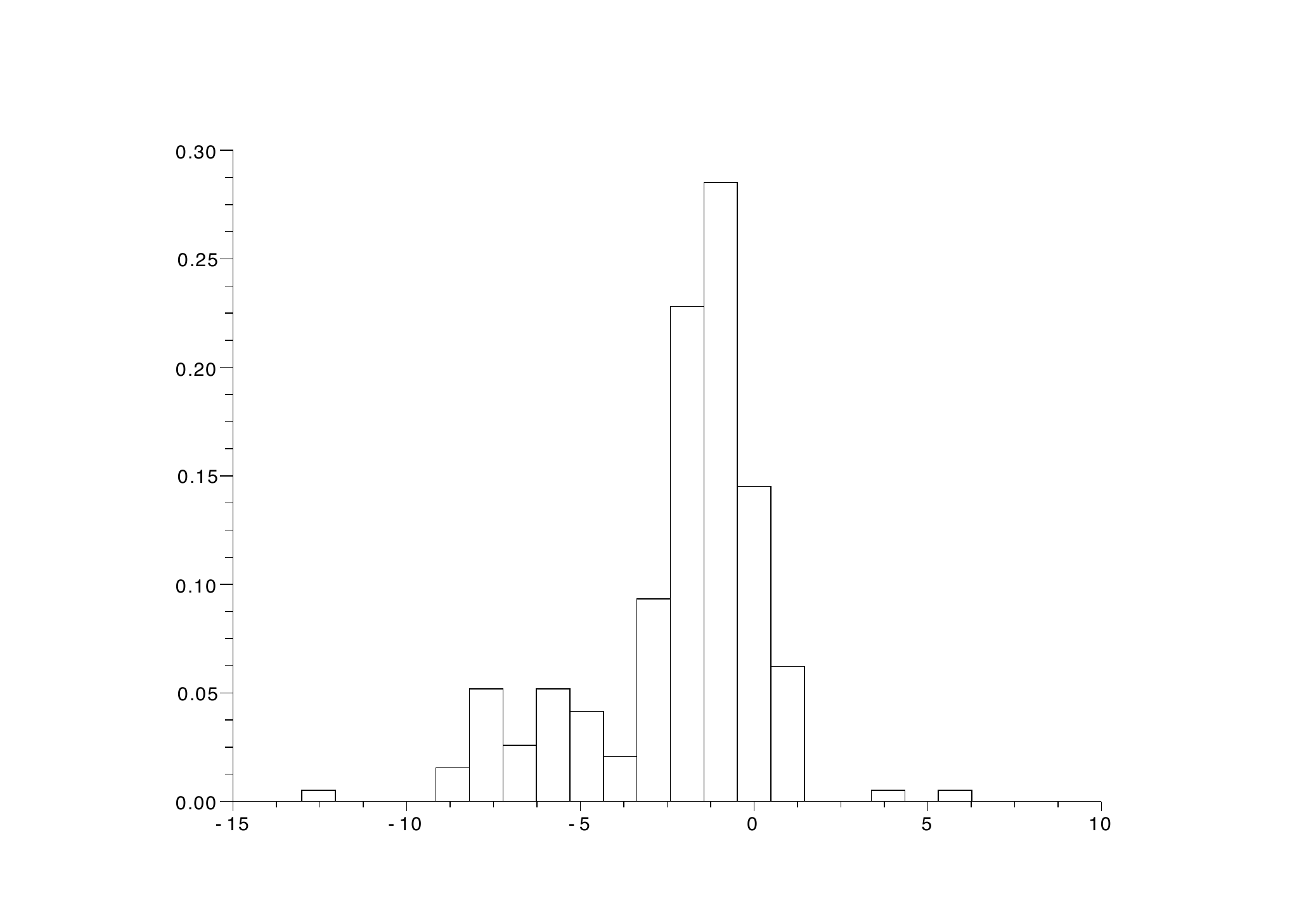}
\caption{ \label{fig:figmixregsim} First row: resp. plot of an original data $(X_i,Y_i)_{1\leq i\leq n}$ according to model (\ref{model}) with  $n=200$,  plot of a wrong $\theta$-tranformation ($\theta\neq \theta^*$), plot of  the true  $\theta^*$-tranformation.  Second row:  resp.  histograms  of the corresponding first row  2nd-coordinate sample data.} 
\end{figure}

We propose in this section  to highlight  the  role played by the  $\theta$-transformation, see (\ref{theta_trans}),  in our  method. For this purpose, we  consider 
an example which corresponds to model (\ref{BDV06model}) taking  $p_*=0.7$,  $\alpha_*=2$, $\beta_*=1$,  $\varepsilon_1^{[j]}\sim {\mathcal N}(0,1),$  $j=0,1$
and $X_1\sim {\mathcal N}(2,3)$. 
In Fig. \ref{fig:figmixregsim} we plot successively a simulated data set  $(X_i,Y_i)_{1\leq i\leq n}$, corresponding to the previous description with $n=200$, 
and the  two   $\theta$-transformed datasets  obtained with $\theta=(1,0.5)$ and $\theta=\theta^*=(2,1)$.  
These figures are completed by adding their corresponding 2nd-coordinate sample data histograms.
Note that these histograms are empirical estimates of the densities $f_\theta$, by  formula (\ref{proj_dens}),  with $\theta$ respectively equal to $(0,0)$, $(1,0.5)$ and $(2,1)$.
We see clearly through these three situations how a progressive transformation of the data allows  one to reach a tractable situation in the sense that it  looks strongly like   the semiparametric contamination  model  (\ref{BDV06model})  studied  in \cite{BDV06} and \cite{BV10} where a known density is mixed with a symmetric unknown density, which corresponds to  the behavior observed in  the second row, third column histogram in Fig.  \ref{fig:figmixregsim}. Loosely speaking the second  idea of our method consists in arguing that once $\theta$ is close to $\theta_*$ we are allowed to estimate the proportion $p$ according to a \cite{BV10} type-method which corresponds to the minimization step  (\ref{estimateur}). In contrast to this technically satisfying idea, the $\theta$-transformation and the choice of the weight distribution $Q$ introduced in (\ref{contrast_prop}) are two sources of serious  difficulties.  In fact when $\beta_*$ is large and the law of the design data has heavy tails with respect to the tails of $f$, then the $\theta_*$ transformation will move  the points  coming from  the $F_0$-population  and located far from the origin,  to extremely distant positions, which implies intuitively  that the  integral type density involved in  (7) should be extremely  heavily tailed. Thus in order to capture the information contained in the  tails of the $\theta$-transformed data set it is important  to weight  sufficiently   the empirical index of symmetry  $H^2(x; p, \tilde
F_{n,\theta}, \hat J_{n,\theta})$  of    expression  (\ref{emp_contrast})  for large values of $x$, which reduces to choosing  an instrumental distribution $Q$ with non-negligible tails with respect to $F_{\theta_*}$.

\subsection{Otimization procedure and simulation study} 
The aim of this section is to  illustrate graphically,  on a two-dimensionnal  examples, the behavior  of the empirical distance $d_n(p, 0,\beta)$
(the parameter $\alpha$ is assumed to be equal to zero) when $p$ and $\beta$ lie close to the true value of the parameter.  For simplicity 
the parameter will still be denoted $\vartheta:=(p,\beta)$, with $\theta:=\beta$  and $d_n(\vartheta):=d_n(p,0, \beta)$.
The interest of this study is 
to understand closely  the influence of the mixing proportion $p$  and the  regression coefficient  $\beta$ on the shape of the  contrast function $d$ (flatness, sharpness, smoothness, etc.). 
Our  models are denoted M1 and M2 and defined according to (\ref{model_arraybis}) as follows\\

\noindent {\bf M1}: $p_*=0.7$, $\beta_*=1$, $V\sim Q={\mathcal N}(0,4^2)$, $\varepsilon^{[j]}\sim  {\mathcal N}(0,1)$, $X\sim {\mathcal N}(0,3^2)$,\\

\noindent {\bf M2}: $p_*=0.3$, $\beta_*=1$, $V\sim Q={\mathcal N}(0,2^2)$, $\varepsilon^{[j]}\sim  {\mathcal N}(0,1)$, $X\sim {\mathcal N}(0,3^2)$,\\

\noindent where $j=0$, 1.\\

\begin{figure}[h!]
\begin{center}
\includegraphics[width=5.9cm]{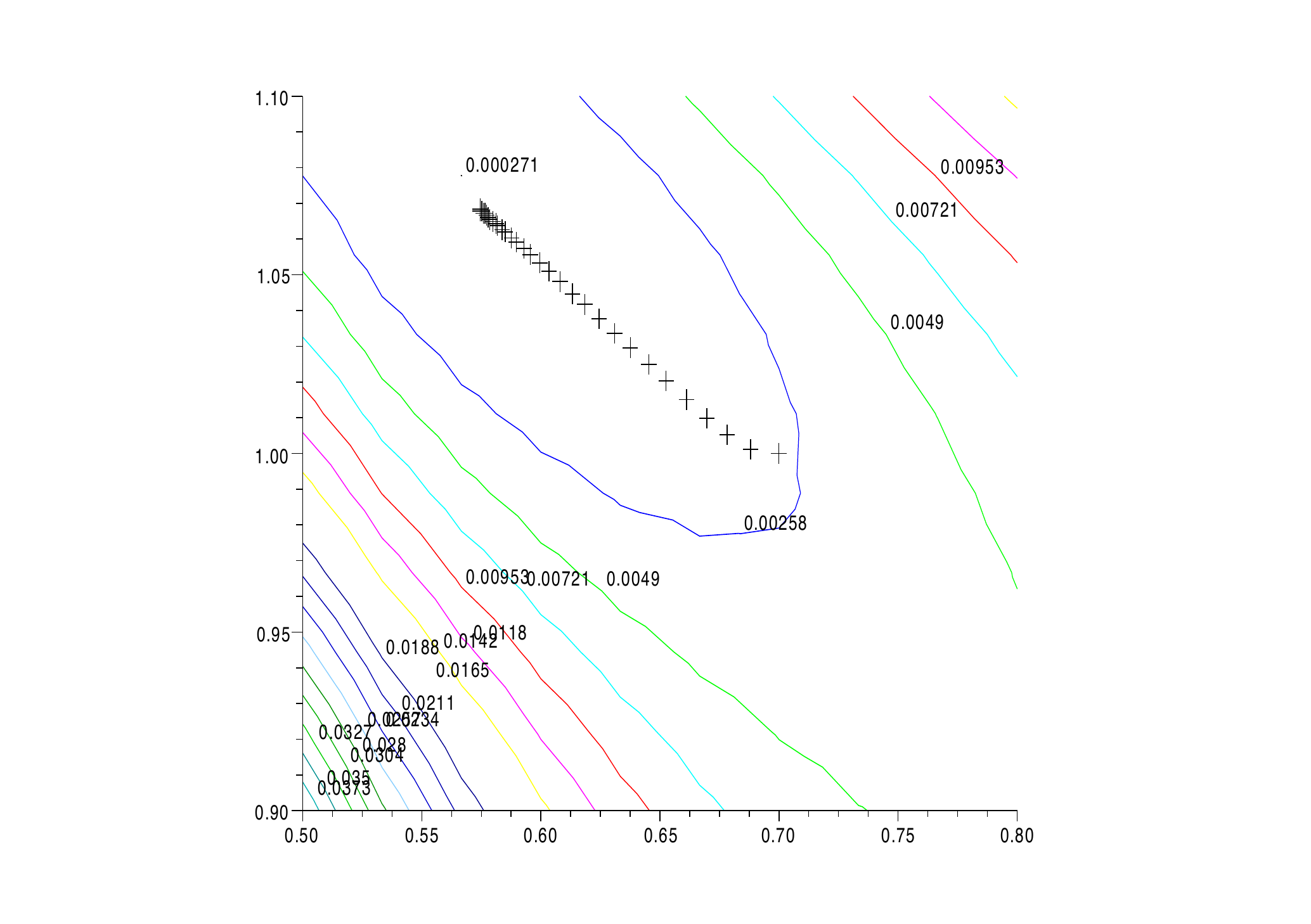}
\includegraphics[width=5.9cm]{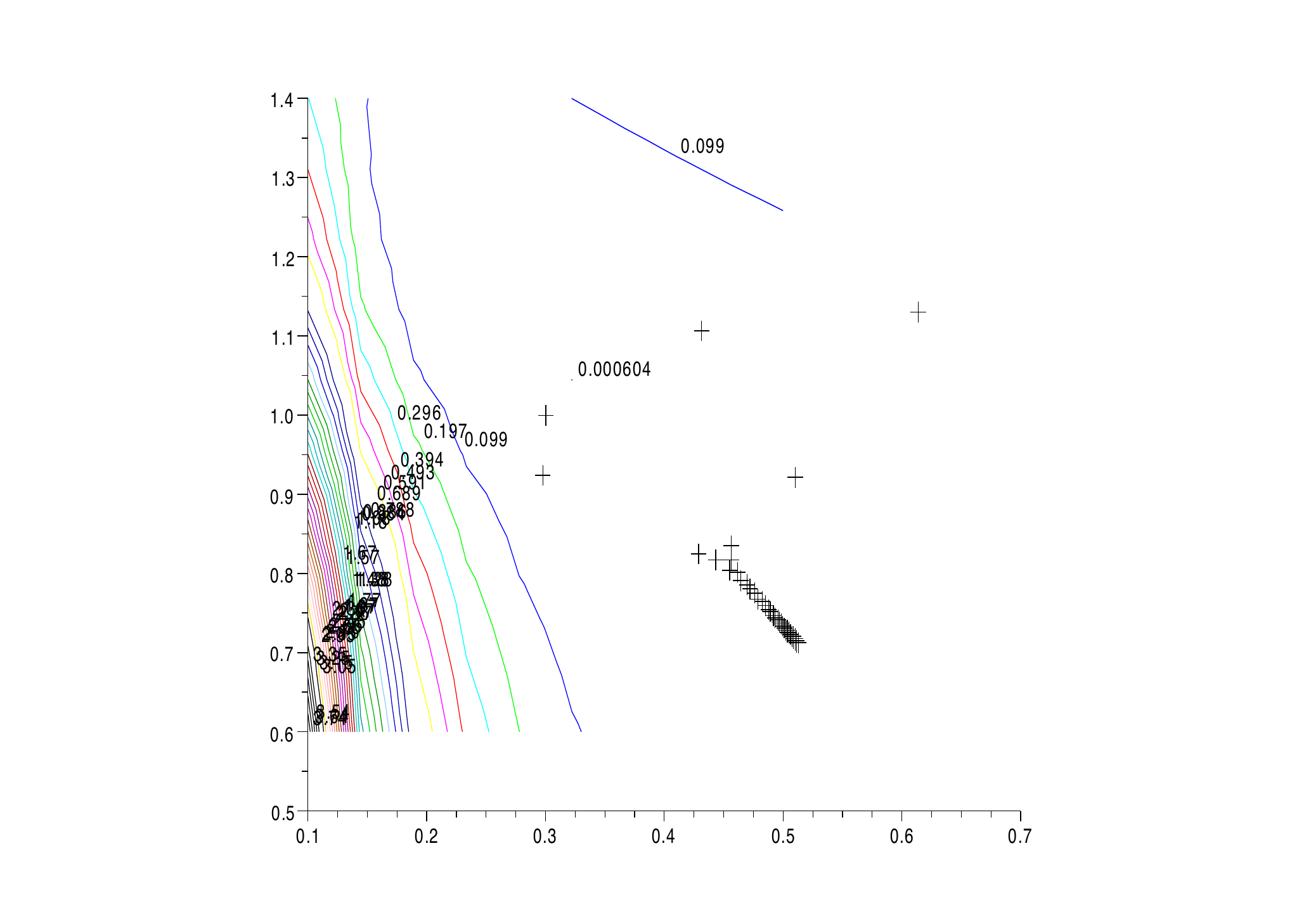}
\end{center}
\caption{\label{fig:fig_shap12} Plot of $(p,\beta)\mapsto d_n(p,0,\beta)$ with  $n=100$,  $\beta_*=1$, $\varepsilon_1^{[j]}\sim {\mathcal N}(0,1),$  $j=0,1$,  $X_1\sim {\mathcal N}(2,3^2)$,
with the difference that  on the left hand side $p_*=0.7$,  $V_1\sim{\mathcal N}(0,4^2)$, when  on the right hand side $p_*=0.3$  $V_1\sim{\mathcal N}(0,2^2)$.} 
\end{figure}

In Fig. \ref{fig:fig_shap12} we plot the mapping $(p,\beta)\mapsto d_n(p,\beta)$ obtained from an M1-sample, resp. M2-sample, of size $n=100$,  where $(p,\beta)\in \Theta_1 =[0.5,0.8]\times[0.9,1.1]$, resp. $(p,\beta)\in \Theta_2= [0.1,0.6]\times[0.6,1.4]$. Notice that according to discussion (CG)  at the end of Section 5.1, model M2 is not necessarily consistently estimated if the parameter space  $\Theta_2$ contains the spurious  solution $\vartheta_{**}=(2p_{*}, \beta_*/2)$, which is voluntary the case here.
In practice, Fig. \ref{fig:fig_shap12} is obtained using  the Scilab  {\texttt{contour2d} } function which plots the  level curves  of $d_n$ evaluated on a  homogeneous $10\times 10$  grid of the rectangular domain $[0.5,0.8]\times[0.9,1.1]$. Fig. \ref{fig:fig_shap12} shows that  the graph of $d_n$ looks like  a sharp  valley with a flat  trough  when $\beta$  is located near $\beta^*$ and $p$ ranges  [0.5,0.8]. Even if on this simulated example  the argmin of $d_n$ is very close to the true value of the parameter, the previous remark suggests   that the estimation of the mixing proportion will be  less robust than the estimation of the regression  coefficient. The observation of the second plot in Fig. \ref{fig:fig_shap12} is more unexpected  since the graph of $d_n$ does not really look like  a contrast function with its high near  $p=0.1$ and its very large and flat trough that covers most of $\Theta_2$ suggesting a  strong lack of robustness of our estimating method in that kind of situation. 
To validate these thoughts we propose to apply a large sample study on the example 
and a third intermediary one  obtained by considering $p_*=0.3$ and $V\sim {\mathcal N}(0,4^2)$. The results of this study will be summarized in Table \ref{T1}.  First we   present the numerical approach used to approximate our M-estimator (\ref{estimateur}). \\

\noindent \textit {Gradient algorithm and tuning parameters.}  
The  gradient  optimization procedure (programmed with Scilab) used to compute our M-estimator $\hat\vartheta_n=(\hat p_n,\hat\beta_n)$ is defined as follows:
\begin{enumerate}[(i)]
\item \texttt{Initialization: } $\vartheta_1=\vartheta_*$,  $\vartheta_2=\vartheta_*+\delta$;
\item \texttt{while} $\|\vartheta_2-\vartheta_1\|_2>\epsilon$ \texttt{do} $\vartheta_1=\vartheta_2$ \texttt{and} $\vartheta_2=\vartheta_1-\gamma^T \dot d_n(\vartheta_1)$;
\item \texttt{else} $\hat \vartheta_n=\vartheta_2$,
\end{enumerate} 
where  $\delta\in \R^2$ is used to create a small perturbation of the initial value, $\epsilon>0$ defines the wanted  stabilization level  in the stopping  algorithm procedure,  and  $\gamma\in {\R^{+*}}^2$  is  a  scale  parameter that needs to be hand-tuned  for good efficiency in practice (to avoid  reverberation phenomena when the score function becomes abruptly sharp). 
The score function   $\dot d_n:=\left(\frac{\partial}{\partial p}d_n, \frac{\partial}{\partial \beta}d_n\right)^T$ 
can  be expressed  into a closed form, {\it i.e.}
\begin{eqnarray*}
\frac{\partial}{\partial p}d_n(\vartheta)&=&  2 \int_\R h_{n,p}(y,\vartheta) H_n (y,\vartheta) dQ_n(y),\\
\frac{\partial}{\partial \beta}d_n(\vartheta)&=&  2\int_\R h_{n,\beta}(y,\vartheta) H_n (y,\vartheta) dQ_n(y) ,
\end{eqnarray*}
where  for all $y\in \R$,
\begin{eqnarray*}
h_{n,p}(y,\vartheta)&:=&\frac{\partial}{\partial p}H_n(y,\vartheta)=-\frac{1}{p^2} \left(  \tilde F_{n,\beta}(y)+ \tilde F_{n,\beta}(-y)-[\hat J_{n,\beta}(y)+\hat J_{n,\beta}(-y)]\right ) \\
h_{n,\beta}(y,\vartheta)&:=&\frac{\partial}{\partial \beta}H_n(y,\vartheta)=\frac{1}{p}\left( \tilde \Psi_{n,\beta}(y)+ \tilde \Psi_{n,\beta}(-y) \right) -\frac{1-p}{p} \left( {j}_{n,\beta}(y)+ {j}_{n,\beta}(-y)\right ) ,
\end{eqnarray*}
and where, from  (\ref{fdr_estim}) and (\ref{dens_estim})
\begin{eqnarray*}
 \tilde \Psi_{n,\beta}(y)&:=&\frac{\partial}{\partial \beta}\tilde F_{n,\beta}(y)\\
 &=&\frac{1}{nb_n}\sum_{i=1}^n\frac{\partial}{\partial \beta}\int_{-\infty}^y K\left(\frac{t-(Y_i-\beta X_i)}{b_n}\right)dt\\
 &=&\frac{1}{n}\sum_{i=1}^n\frac{\partial}{\partial \beta}\int_{-\infty}^{\frac{y+\beta X_i-Y_i}{b_n}} K\left(u\right)du\\
 &=&\frac{1}{nb_n}\sum_{i=1}^nX_iK\left(\frac{y+\beta X_i-Y_i}{b_n}\right),
 \end{eqnarray*}
 and similarly, from  (\ref{J-estim-final})
\begin{eqnarray*}
 {j}_{n,\beta}(y)&:=&\frac{\partial}{\partial \beta}\hat J_{n,\beta}(y)\\
 &=&\frac{1}{n}\sum_{i=1}^nX_i  f_0\left(y+\beta X_i\right).
\end{eqnarray*}
The kernel $K$ used to compute (\ref{dens_estim}), is a triangular kernel defined by
$$
K(x)=(1-|x|)\un_{-1\leq x\leq 1},\quad\quad x\in \R,
$$
and the bandwidth $b_n=\sqrt{1+4p(1-p)}(4/(3n))^{1/5}$ (proposed by \cite{BA03} for gaussian distributions and implemented  in  \texttt{R}), both obviously satisfying condition (K). The results summarized in Table 1 were obtained with the following hand-tuned parameters:  $\delta=0.01$, $\epsilon=0.005$,  $\gamma=(0.2,0.5)^T$, and an example of  stabilization for this set
of tuning parameters is illustrated in Fig  \ref{fig:fig_shap12}, where  the successive positions (until stabilization) of our algorithm  are depicted by  cross symbols.  \\

\begin{table}[h!]
\caption{\label{T1}Mean and Std. Dev. of 100 estimates of $p$, $\beta$. }
\centerline{\begin{tabular}{cccc}
\hline
$n$  & $(p_*,\beta_*,\sigma_V)$  & Empirical means & Standard deviation \\
\hline
100 &  (0.7,1,4) & (0.7055,1.0051) & (0.0373,0.0697) \\
200 &  (0.7,1,4) & (0.6976,0.9965) & (0.0307,0.0590) \\
500 &  (0.7,1,4) & (0.6954,1.0059) & (0.0296,0.0358) \\
100 &  (0.3,1,4) & (0.3100,0.9581) & (0.0577,0.1252) \\
200 &  (0.3,1,4) & (0.2965,0.9851) & (0.0501,0.0855 ) \\
500 &  (0.3,1,4) & (0.2975,1.0178) & (0.0284,0.0414) \\
100&  (0.3,1,2)  & (0.3971, 0.8587) & (0.0942, 0.2213)\\
200&  (0.3,1,2)  & (0.3982,0.9149) & (0.0835,0.1900)\\
500& (0.3,1,2)  & (0.3315, 0.9683) & (0.0524, 0.1067)\\
\hline
\end{tabular}}
\end{table}

\noindent \textit{Comments on Table \ref{T1}}.
First of all, it is interesting to compare the performances summarized in rows 1--3  of  Table 1 to those obtained in \cite{BV10}, p. 35, Table 1 where the 
model of interest is  (\ref{BDV06model}),  with $p=0.7$, $\mu=3$, and $f_0$ and $f$ are respectively the pdfs corresponding to the $\mathcal{N}(0,1)$ and $\mathcal{N}(0,(1/2)^2)$ distributions. Even if these two  models are not  strictly comparable we think that it is interesting, in order to  highlight the drawbacks induced by  the $\theta$-transformation and the choice of $Q$ discussed above, to compare pairwise the performance obtained on the  mixing proportion $p$  and  the parameters that influence the location of the $F$-population, \textit{i.e.}  $\beta$ and $\mu$. From the numerical point of view, we  easily check that the bias of our estimators,  for  both models,  is negligible. However  it also appears that the standard deviation associated
to $(\hat p_n, \hat \beta_n)$ decreases significantly slower than  the standard deviation associated to $(\hat p_n, \hat \mu_n)$ when $n$ grows.  The performance summarized in rows 4--6 of Table 1, which corresponds to $p=0.3$ (and hence signifies that the population that will move far from its original position due to the  $\theta$-transformation will be more important), is 
very instructive. We observe that for small  $n$ ($n=100$, 200) the standard deviations associated to $(\hat p_n, \hat \beta_n)$ are  dramatically large compared to those obtained
with $p=0.7$.  Let $\mbox{std}(n, p_*,\beta_*,\sigma_V)$ the couple of standard deviations calculated in the last column of Table 1 under $(n, p_*,\beta_*,\sigma_V)$.
If we compute componentwise    the  ratios  $\mbox{std}(n,0.3,1,4)/\mbox{std}(n,0.7,1,4)$    respectively for $n=100, 200, 500$ we obtain approximately $(1.54,1.8)$, $(1.67, 1.44)$, and $(0.95, 1.17)$ which seems to suggest  that when $n$ becomes large 
the side effect of the $\theta$-transformation vanishes (probably thanks to the size of $n$, which increases globally the precision of the empirical estimates, and the tails of $Q$, that allow the algorithm to take these improvements into account efficiently).
The performances summarized in rows 7--9 of Table 1, seems to confirm the concerns expressed about model M2. We recall  that model M2 is badly affected by the two following drawbacks : smallness of $p_*$ (synonymous with  important population shifted far  by the $\theta$-transformation and  existence of a spurious  solution) and a smallness of $\sigma_V$ which is then clearly not sufficient to counteract the smallness of $p_*$ (and its consequences).
We think in particular that, in model M2, the empirical contrast $d_n$ is more easily closer to 0 under $\beta_{**}=\beta_*/2$ since as  explained in Section 4.1., this value is then significantly smaller than $\beta_*$. This last remark explains why, in spite of the fact that our algorithms were initialized at the true  parameter value, our estimates are strongly biased (attracted quite often by the spurious solution $\vartheta_{**}$).\\

\begin{figure}[h!]
\begin{center}
\includegraphics[width=3cm]{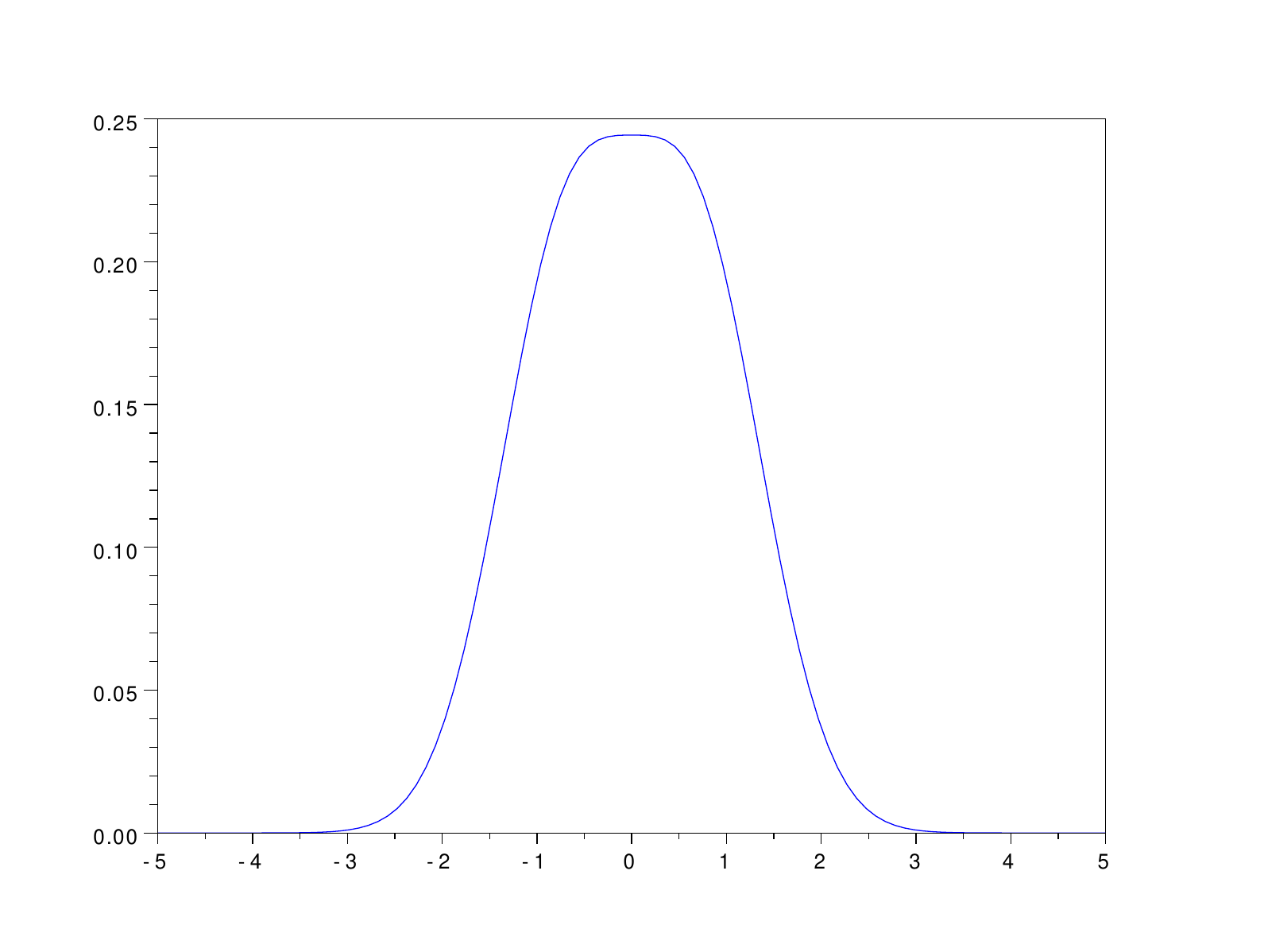}
\includegraphics[width=3.2cm]{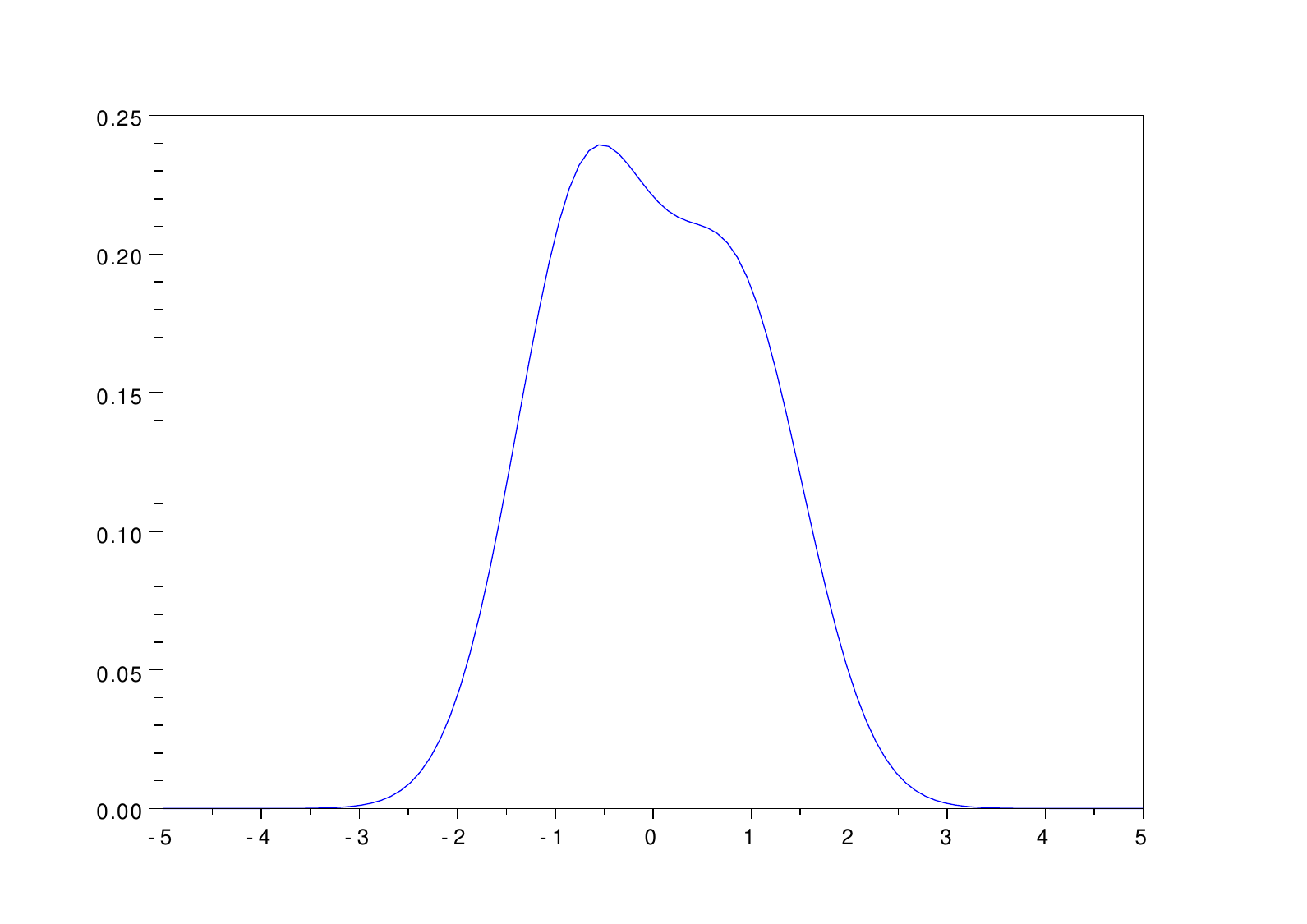}
\includegraphics[width=3.3cm]{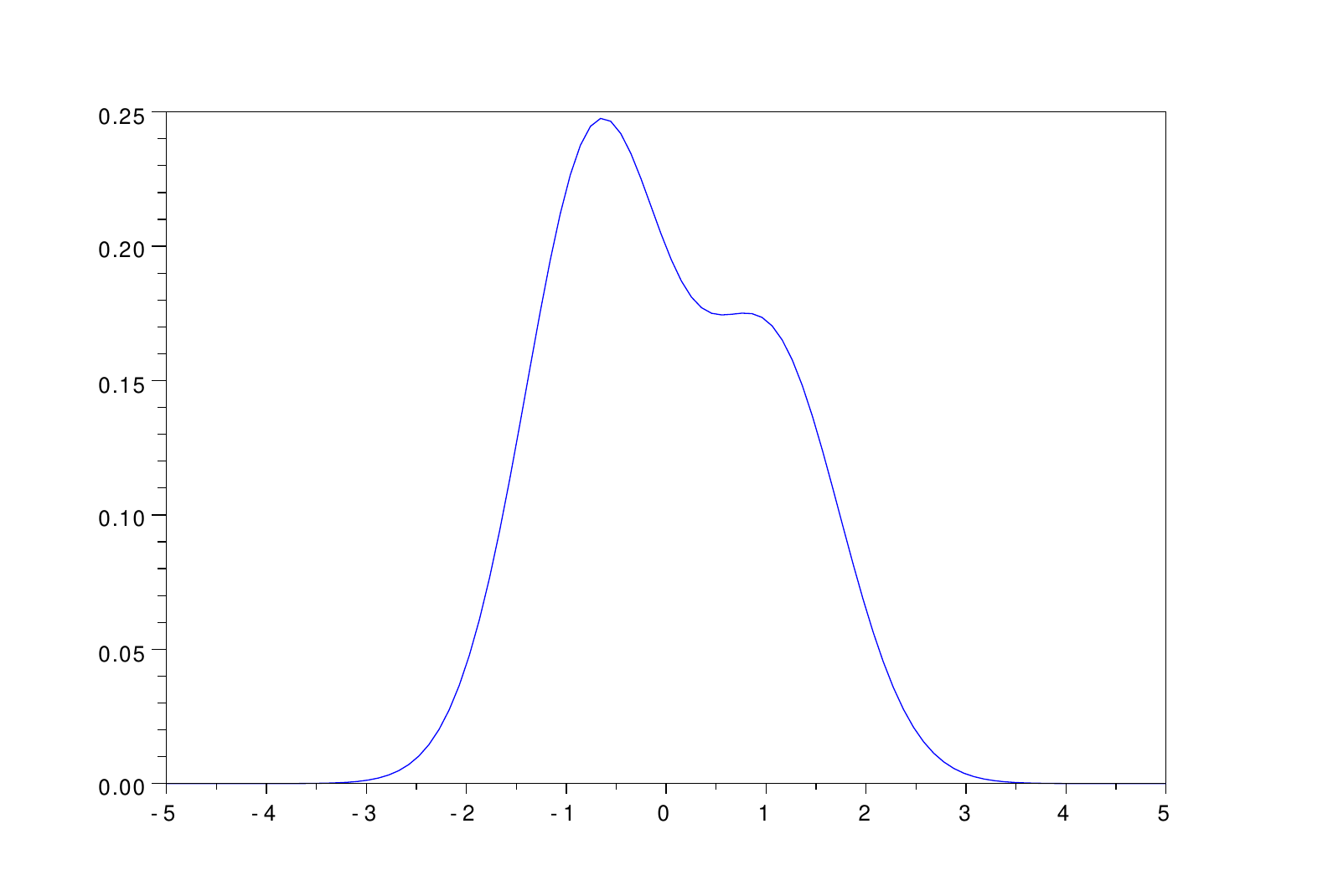}
\end{center}
\caption{\label{fig:figgraphnonsym} Graphs of  the pdfs  corresponding to  the  mixture distribution  $\lambda {\mathcal N}(-0.7,1/\sqrt{2})+(1-\lambda) {\mathcal N}(0.7\lambda/(1-\lambda),1/\sqrt{2})$, obtained by considering $\lambda=0.5, 0.55, 0.6$.} 
\end{figure}

\noindent\textit{Robustness with respect to the symmetry assumption}. We propose to conclude this simulation study  by  testing our method in situations where the law of the error $\varepsilon^{[1]}$ is no longer symmetric. For this purpose we consider  again model M1 and replace the distribution of $\varepsilon^{[1]}$ by the  mixture $$\lambda {\mathcal N}(-0.7,1/\sqrt{2})+(1-\lambda) {\mathcal N}(0.7\lambda/(1-\lambda),1/\sqrt{2})$$ which pdf, denoted $f_\lambda$,  is nonsymmetric if $\lambda\neq 0.5$ but has  a mean equal to 0 and a variance equal to 0.5
for all $\lambda\in(0,1)$. In our simulations we consider   successively  $\lambda=0.5, 0.55, 0.6$  which leads to consider
pdfs for $\varepsilon^{[1]}$  which  graphs are plotted in Fig \ref{fig:figgraphnonsym}. 
Some performances of our method on these examples are summarized in Table \ref{T2}.\\
\begin{table}[h!]
\caption{\label{T2}Mean and Std. Dev. of 100 estimates of $p$, $\beta$. }
\centerline{\begin{tabular}{cccc}
\hline
$n$  & $\lambda$  & Empirical means & Standard deviation \\
\hline
100 & 0.5  & (0.7035,1.0229) & (0.0427,0.0814) \\
200 & 0.5  & (0.7012,1.0068) & (0.0390,0.0774) \\
500 & 0.5 & (0.6997,1.0059) & (0.0244,0.0488) \\
100 & 0.55  & (0.6854,1.0837) & (0.0485,0.0858) \\
200 & 0.55 & (0.6890,1.0805) & (0.0431,0.0716) \\
500 & 0.55 & (0.6922,1.0699) & (0.0377,0.0519) \\
100& 0.6 &  (0.6731,1.1314) & (0.0543,0.0952)\\
200& 0.6 &  (0.6693,1.1061) & (0.0490,0.0868)\\
500& 0.6 & (0.6775,1.0928) & (0.0392,0.0557)\\
\hline
\end{tabular}}
\end{table}

\noindent \textit{Comments on Table \ref{T2}}. Note that when  $\lambda=0.5$  (symmetric case)   the performances of our method  are very close to those obtain on model M1. However  for $n=100, 200$ the standard deviation of our estimates is larger than those obtained in the M1 model, when for $n=500$ the standard deviation becomes slightly smaller.  This behavior can probably be explained by the fact  that the graph of $f_{0.5}$ is flat on its top which intuitively do not help much in locating  the  axis of symmetry for small values of $n$.   On the other hand we can expect that for  $n=500$, helped with the fact that $\mbox{var}(\varepsilon^{[1]})$ is here equal to $0.5$  when it was equal to $1$ in M1,  our nonparametric estimators perform better than in model M1  which should explain the good performances observe in the  third row of Table \ref{T2}. When $\lambda=0.55, 0.6$ it appears that the parameter $\beta$ is always overestimated. This phenomenon can be explained by the fact that our method try to determine a pseudo-axis of symmetry  adapted to the  shapeless graph of  $f_\lambda$  which qualitatively is placed on the left side of the origin. This remark implies that  the $\theta$-transformation needed to transform the first integral in (\ref{proj_dens}) into an almost even density (see Fig. \ref{fig:fig_shap12}) have to contain a $\beta$ greater than $\beta_*$.

\section{Appendix}

\subsection{\it  Conditions (R) and (C) in the  Gaussian Case}\label{ident_gauss_sect}
In this section we discuss conditions (R) and (C) when the true underlying model is  a contaminated
Gaussian  regression model with Gaussian design, {\it i.e.}, $f$, $f_0$, $h$ are respectively the pdfs 
of the ${\mathcal N}(0,m)$, ${\mathcal N}(0,m_0)$, and ${\mathcal N}(E(X),\sigma_h^2)$ distributions. \\

\noindent\textit{Comments on  condition (R).} Conditions (R) i-iii)  are standard and easy to verify in the above model.
On the other hand,  it is interesting to show how conditions (R) iv-v) arise naturally in this case. \\

\noindent{\it Condition (R) iv).} We  show for simplicity  that the first condition in (R) iv) (the same kind of proof works  also  for the second one)  holds  when
$i=0$, $\theta^*=(\alpha_*,\beta_*)\in {\R^{+*}}^2$ and $m_0=1$. We write the decomposition
\begin{eqnarray*}
&&|F_0(y+\theta_*\odot x)-F_0(y-\theta_*\odot x)|=\\
&&~~~~~~~~~~~~~~~~~|F_0(y+\theta_*\odot x)-F_0(y-\theta_*\odot x)|\un_{y>1-\theta_*\odot x, ~x<-\frac{\alpha_*}{\beta_*}}\\
&&~~~~~~~~~~~~~+|F_0(y+\theta_*\odot x)-F_0(y-\theta_*\odot x)|\un_{y>1+\theta_*\odot x, ~x\geq-\frac{\alpha_*}{\beta_*}}\\
&&~~~~~~~~~~~~~+|F_0(y+\theta_*\odot x)-F_0(y-\theta_*\odot x)|\un_{y<-1+\theta_*\odot x, ~x<-\frac{\alpha_*}{\beta_*}}\\
&&~~~~~~~~~~~~~+|F_0(y+\theta_*\odot x)-F_0(y-\theta_*\odot x)|\un_{y<-1-\theta_*\odot x, ~x\geq -\frac{\alpha_*}{\beta_*}}\\
&&~~~~~~~~~~~~~+|F_0(y+\theta_*\odot x)-F_0(y-\theta_*\odot x)|\un_{-1+\theta_*\odot x\leq y\leq 1-\theta_*\odot x, ~x<-\frac{\alpha_*}{\beta_*}}\\
&&~~~~~~~~~~~~~+|F_0(y+\theta_*\odot x)-F_0(y-\theta_*\odot x)|\un_{-1-\theta_*\odot x\leq y\leq 1+\theta_*\odot x, ~x\geq -\frac{\alpha_*}{\beta_*}}.
\end{eqnarray*}
Consider the first term on  the right hand side of the above decomposition (the three following terms being treated in entirely same way). For all $y>1-\theta_*\odot x$ with  $x<-\alpha_*/\beta_*$ we  have
$y-\theta_*\odot x>y+\theta_*\odot x>1$.  Since  for $t>0$ large enough, the  inequality (\ref{majcdfGauss})  is valid 
\begin{eqnarray}\label{majcdfGauss}
\frac{\exp(-t^2) }{\sqrt{2\pi}} \left(\frac{1}{t} -\frac{1}{t^3} \right) \leq 1-F_0(t)\leq \frac{\exp(-t^2) }{t\sqrt{2\pi}},
\end{eqnarray}
we have in particular that for all  $t>1$, $0 \leq 1-F_0(t)\leq \exp(-t^2)/\sqrt{2\pi}$.
Hence it follows  that for $y>1-\theta_*\odot x$ with  $x<-\alpha_*/\beta_*$:
\begin{eqnarray*}
|F_0(y+\theta_*\odot x)-F_0(y-\theta_*\odot x)|\leq  \frac{\exp(-(y+\theta_*\odot x)^2)+ \exp(-(y-\theta_*\odot x)^2)}{\sqrt{2\pi}},
\end{eqnarray*}
which proves that this first term is $h(x)dxdy$ integrable.
Let us now sum  the  last two terms of the above decomposition and notice that
\begin{eqnarray*}
&&|F_0(y+\theta_*\odot x)-F_0(y-\theta_*\odot x)|\\
&&\times \left(\un_{-1+\theta_*\odot x\leq y\leq 1-\theta_*\odot x, ~x<-\frac{\alpha}{\beta_*}}+\un_{-1-\theta_*\odot x\leq y\leq 1+\theta_*\odot x, ~x\geq -\frac{\alpha}{\beta_*}}\right)\\
&\leq& 2 \un_{-1-|\theta_*\odot x|\leq y\leq 1+|\theta_*\odot x|}.
\end{eqnarray*}
We thus  prove that this sum of terms  is also  $h(x)dxdy$ integrable.\\

\noindent{\it Condition (R) v).} We consider for simplicity
the construction of the bounding function $\ell^0_{1,1}$ when $(\alpha,\beta)\in \Phi= [\underline{\alpha}, \overline{\alpha}]\times [\underline{\beta}, \overline{\beta}]$, with $(\underline{\alpha}, \underline{\beta})\in {\R^{+*}}^2$ and $m=m_0=1$.
Notice first that  for all $(x,y)\in \R^2$
$$
|f_0^{(1)}(y+\theta\odot x)|\leq \frac{|y|+\overline{\alpha}+\overline{\beta}|x|}{\sqrt{2\pi}}\exp\left(-\frac{(y+\alpha+\beta x)^2}{2}\right).
$$
Secondly it is easy to check that  for all $(x,y)\in \R^2$ and all $\theta\in \Phi$:
\begin{eqnarray*}
&&\exp\left(-\frac{(y+\alpha+\beta x)^2}{2}\right)\\
&&~~~~~~~~~~~~~~~~~\leq \exp\left(-\frac{(y+\underline{\beta}x)^2}{2}\right)\un_{x\geq 0, y\geq 0}\\
&&~~~~~~~~~~~~~~~~~~ + \exp\left(-\frac{\min(|y+\underline\alpha+\underline{\beta}x|, |y+\overline\alpha+\overline{\beta}x|)^2}{2}\right)\un_{x\geq 0, y\leq 0} \\
&&~~~~~~~~~~~~~~~~~~+\exp\left(-\frac{\min(|y+\underline\alpha+\overline{\beta}x|, |y+\overline\alpha+\underline{\beta}x|)^2}{2}\right)\un_{x\leq 0, y\in \R}\\
&&~~~~~~~~~~~~~~~~~\leq B_{\Phi}(x,y),
\end{eqnarray*}
where 
\begin{eqnarray*}
B_{\Phi}(x,y)&:=&
\exp\left(-\frac{(y+\underline{\beta}x)^2}{2}\right)
+\exp\left(-\frac{(y+\underline\alpha+\underline{\beta}x)^2}{2}\right)\\
&&+ \exp\left(-\frac{(y+\overline\alpha+\overline{\beta}x)^2}{2}\right)
+ \exp\left(-\frac{(y+\underline\alpha+\overline{\beta}x)^2}{2}\right)\\
&&+ \exp\left(-\frac{(y+\overline\alpha+\underline{\beta}x)^2}{2}\right).
\end{eqnarray*}
Thus we can propose 
$\ell^0_{1,1}(x,y)=|x|(|y|+\overline{\alpha}+\overline{\beta}|x|)/\sqrt{2\pi}B_\Phi(x,y) \exp(-x^2/2)$,  which  clearly belongs to  $L_1(\R^2)$, as a candidate for the uniformly bounding function 
satisfying condition (R) v).\\

\noindent\textit{Comments on condition (C).} In the whole Gaussian case,  expression of (\ref{Fourier_eq}) becomes:
\begin{eqnarray}\label{ident_gauss}
&& p_*\sin((\alpha-\alpha_*)+(\beta-\beta_*)E(X))t) \exp\left(-\frac{t^2}{2}(\sigma_h^2(\beta-\beta_*)^2+m)\right)\nonumber\\
&&~~~~~~~~~~=(p_*-p)\sin((\alpha+\beta E(X))t)\exp\left(-\frac{t^2}{2}(\sigma_h^2\beta^2+m_0)\right).
\end{eqnarray}
We suppose first that $p\neq p_*$, and denote $\xi:=\alpha+\beta E(X)$, $\xi_*:=\alpha_*+\beta_* E(X)$, $\Sigma_{\beta-\beta_*}:=\sigma_h^2(\beta-\beta_*)^2+m$,  and  $\Sigma_{\beta}:=\sigma_h^2\beta^2+m_0$. 
Taking the first and third  order derivative of  (\ref{ident_gauss}) at point $t=0$ we get
the  conditions
 \begin{eqnarray}\label{ident_gauss0}
p_*\xi_*=p\xi,\quad\mbox{and}\quad (p-p_*)\xi(3\Sigma_{\beta-\beta_*}+\xi^2)+p_*(\xi-\xi_*)(3 \Sigma_{\beta}+(\xi-\xi_*)^2)=0.
\end{eqnarray} 
Introducing  the first relation in (\ref{ident_gauss0}) into  the second one,  we obtain
\begin{eqnarray}\label{ident_gauss1}
 \frac{p-p_*}{p}\xi_*p_*\left(3[\Sigma_{\beta}-\Sigma_{\beta-\beta_*}]+\frac{2p_*p-p^2}{p^2}\xi_*^2\right)=0.
\end{eqnarray}
Now we observe that, to insure the validity of expression (\ref{ident_gauss}), the factors multiplied by   the $\sin$ terms on both sides of (\ref{ident_gauss}) must be, at least,  equivalent as $t\rightarrow \infty$.
This last remark implies that $\Sigma_{\beta}=\Sigma_{\beta-\beta_*}$,  or equivalently $\beta=\beta_*/2+ (m_0-m)/2\beta_*\sigma_h^2$, and thus  (\ref{ident_gauss1}) leads to
\begin{eqnarray}\label{ident_gauss2}
 p=2p_*.
\end{eqnarray}
Using now the first relation in (\ref{ident_gauss}) and (\ref{ident_gauss2}), we then obtain $\alpha=\alpha_*/2+E(X)(m_0-m)/2\beta_*\sigma_h^2$.\\

\noindent The consequences of the previous comments  can be presented as follows:\\

\noindent \textit{Discussion (CG)}:

\begin{enumerate}[i)]
\item If $p_*>1/2$ then the set of parameters $\vartheta\in \Theta$ satisfying condition (\ref{ident_gauss}) is always empty, since  $p=2p^*>1$ 
is not an admissible solution.
\item If $p_*\leq 1/2$ and if, for example,  $E(X)=0$ and $m_0=m$ then $\vartheta_{**}=(2p^*,\alpha_*/2, \beta_*/2)$. In such 
a case it would be crucial  to build a conveniently constrained  parametric space (most of the time  a plot of the dataset  
 helps in  building   reasonnable constraints on the intercept and slope parameter spaces) expecting that it contains $\vartheta_*$ but not $\vartheta_{**}$. 
 \item More generaly one can expect that when the shape of the sample data (see {\it e.g.} Fig. \ref{fig:figmixregsim})  suggest  that
$(m_0-m)/2\beta_*\sigma_h^2$ is negligible with respect to $\alpha_*$ and $\beta_*$, which occurs when $m_0$ is close to $m$ or/and $\sigma_h^2\beta_*$ is very large,
then the solution proposed in ii) is loosely speaking still valid. 
\end{enumerate}

\subsection{\it Explicit formula of $H(\cdot,\vartheta)$ and its derivatives} \label{expl_H}
In this section all the expressions are valid for all $(\vartheta,y)\in \Theta\times\R$, and the computation of the various
derivative functions (under the integral sign) are all allowed according to  Lebesgue's Theorem and condition (R).
According to (\ref{proj_dens}) and (\ref{inv_formula}), we have 
\begin{eqnarray*}
H(y,\vartheta)&=&\frac{p^*}{p}\left( \int_{-\infty}^y\int_\R f(z+(\theta-\theta^*)\odot x) h(x) dxdz\right.\\
&&+\left. \int_{-\infty}^{-y}\int_\R f(z+(\theta-\theta^*)\odot x) h(x) dxdz\right)\\
&&+\frac{p-p^*}{p}\left( \int_{-\infty}^y\int_\R f_0(z+\theta\odot x) h(x) dxdz\right.\\
&&+\left.\int_{-\infty}^{-y} \int_\R f_0(z+\theta\odot x) h(x) dxdz\right)-1.
\end{eqnarray*}
For simplicity we introduce 
\begin{eqnarray*}
F^{\theta}(y)&=& \int_{-\infty}^y\int_\R f(z+(\theta-\theta^*)\odot x) h(x) dxdz,\\
F_0^{\theta}(y)&=& \int_{-\infty}^y\int_\R f_0(z+\theta\odot x) h(x) dxdz,
\end{eqnarray*}
which leads to
\begin{eqnarray*}
\frac{\partial}{\partial p}H(y,\vartheta)=-\frac{p^*}{p^2}\left[ (F^\theta(y)+F^\theta(-y))-(F_0^\theta(y)+F_0\theta(-y))\right].
\end{eqnarray*}
Let us denote 
\begin{eqnarray*}
\dot F^\alpha(y)&:=&\frac{\partial}{\partial \alpha}F^\theta(y)=\int_{-\infty}^y \int_\R \dot f(z+(\theta-\theta_*)\odot x)h(x)dxdz,\\
\dot F_0^\alpha(y)&:=&\frac{\partial}{\partial \alpha}F_0^\theta(y)=\int_{-\infty}^y \int_\R \dot f_0(z+\theta\odot x)h(x)dxdz,
\end{eqnarray*}
and for
\begin{eqnarray*}
\dot F^\beta(y)&:=&\frac{\partial}{\partial \beta}F\theta(y)=\int_{-\infty}^y \int_\R x \dot f(z+(\theta-\theta^*)\odot x)h(x)dxdz,\\
\dot F_0^\beta(y)&:=&\frac{\partial}{\partial \beta}F_0\theta(y)=\int_{-\infty}^y \int_\R x\dot f_0(z+\theta\odot x)h(x)dxdz,
\end{eqnarray*}
we obtain
\begin{eqnarray*}
\frac{\partial}{\partial \alpha}H(y,\vartheta)=\frac{p^*}{p}\left(\dot F^\alpha(y)+\dot F^\alpha(-y) \right)+\frac{p-p^*}{p}\left(\dot F_0^\alpha(y)+\dot F_0^\alpha(-y) \right).
\end{eqnarray*}

\begin{eqnarray*}
\frac{\partial}{\partial \beta}H(y,\vartheta)=\frac{p^*}{p}\left(\dot F^\beta(y)+\dot F^\beta(-y) \right)+\frac{p-p^*}{p}\left(\dot F_0^\beta(y)+\dot F_0^\beta(-y) \right).
\end{eqnarray*}
At point $\vartheta=\vartheta_*$  the Hessian matrix of $H(\cdot,\vartheta)$  defined in (\ref{hess_theta*}) is obtained  by considering 
\begin{eqnarray}\label{deriv H}
\dot H(y,\vartheta_*)=\left(
\begin{array}{c}
 \displaystyle{\frac{1}{p}\left(F_0^{\theta_*}(y)+F_0^{\theta_*}(-y)-1\right)}\\
2f(y)\\
2f(y)E(X)
\end{array}
\right).
\end{eqnarray}
Let us denote now
\begin{eqnarray*}
\ddot F^{\beta,\alpha}(y)=\frac{\partial^2 }{\partial \beta\partial\alpha } F^\theta(y)&=&\int_{-\infty}^y \int_\R x\ddot f(z+(\theta-\theta_*)\odot x)h(x)dxdz,\\
\ddot F_0^{\beta,\alpha}(y)=\frac{\partial^2 }{\partial \beta\partial \alpha } F_0^\theta(y)&=&\int_{-\infty}^y \int_\R x \ddot f_0(z+\theta\odot x)h(x)dxdz,\\
\ddot F^{\alpha,\alpha}(y)=\frac{\partial^2 }{\partial \alpha^2 } F^\theta(y)&=&\int_{-\infty}^y \int_\R \ddot f(z+(\theta-\theta_*)\odot x)h(x)dxdz,\\
\ddot F_0^{\beta,\beta}(y)=\frac{\partial^2 }{\partial \beta^2 } F_0^\theta(y)&=&\int_{-\infty}^y \int_\R x^2 \ddot f_0(z+\theta\odot x)h(x)dxdz.
\end{eqnarray*}
We then obtain 
\begin{eqnarray*}
\frac{\partial^2}{\partial p^2}H(y,\vartheta)&=&\frac{p^*}{p^3}\left[ (F^\theta(y)+F^\theta(-y))-(F_0^\theta(y)+F_0\theta(-y))\right],\\
\frac{\partial^2}{\partial u\partial p}H(y,\vartheta)&=&-\frac{p^*}{p^2}\left[ (\dot F^{u}(y)+\dot F^u(-y))-(\dot F_0^u(y)+\dot F_0^u(-y))\right],\quad u=\alpha,~\beta,\\
\frac{\partial^2}{\partial u\partial v}H(y,\vartheta)&=&\frac{p^*}{p}\left(\ddot F^{u,v}(y)+\ddot F^{u,v}(-y) \right)\\
&&+\frac{p-p^*}{p}\left(\ddot F_0^{u,v}(y)+\ddot F_0^{u,v}(-y) \right),\quad u=\alpha, ~\beta, \quad v=\alpha, ~\beta.
\end{eqnarray*}

\subsection{\it Boundedness}\label{bound}

\noindent{\it Boundedness of $\Psi_\theta (\cdot)$ and $\Psi'_\theta(\cdot)$}. If $f$ and $f_0$  are supposed to be bounded over $\R$ then  we clearly have from  (\ref{proj_dens}) that
$$|\Psi_\theta (y)|\leq \|f\|_\infty+\|f_0\|_\infty,\quad\quad (\theta,y)\in \Phi\times \R.$$
The same kind of argument holds to prove boundedness of $\dot f_\theta(\cdot)$ when (R) ii) is supposed.\\

\noindent {\it Boundedness of $H(\cdot,\vartheta)$}. Since for all $\theta\in \Phi$ the functions $F_\theta(\cdot)$ and $J_\theta(\cdot)$ are both  cdfs,  we thus have, since $\delta \leq p\leq 1-\delta$, from expression (\ref{H_def}):
\begin{eqnarray}\label{boundH}
H(y,\vartheta)\leq \frac{4}{\delta}+1, \quad\quad (y,\vartheta)\in \R\times \Phi.
\end{eqnarray}

\subsection{\it Integrable Lipschitz property  of  $\Psi_\theta (\cdot)$}\label{integ_lipf}
From (\ref{proj_dens}), for  all $(y,\theta,\theta')\in \R\times  \Phi^2$ we have  
\begin{eqnarray}\label{Lip_decomp}
|\Psi_\theta(y)-f\Psi_{\theta'}(y)|&\leq&  p_*\int_{\R} |f(y+(\theta-\theta_*)\odot x)-f(y+(\theta'-\theta_*)\odot x)|h(x)dx\nonumber\\
&&+(1-p_*) \int_{\R} |f_0(y+\theta\odot x)-f_0(y+\theta'\odot x)|h(x)dx.
\end{eqnarray}
Consider for simplicity  the first integral term on the right hand side of (\ref{Lip_decomp})  (the same argument holding for the second term). According to  the
Mean Value Theorem  there exists,  for all $(x,y)\in \R^2$ and $(\theta,\theta')\in \Phi^2$,  a value $\gamma:=\gamma(x,y,\theta,\theta')$  belonging to the line segment with extremities
$y+(\theta-\theta_*)\odot x$ and $y+(\theta'-\theta_*)\odot x$, or equivalently a $\bar\theta:=\bar\theta(x,y,\theta,\theta')$  belonging to the line segment with extremities $\theta$ and $ \theta'$ such that
$\gamma=y+(\bar \theta-\theta_*)\odot x$ and 
\begin{eqnarray*}
 &&|f(y+(\theta-\theta_*)\odot x)-f(y+(\theta'-\theta_*)\odot x) |\\
 &&~~~~~~~~~~~~~~=|\dot f(\gamma) (\alpha-\alpha'+(\beta-\beta')x)|\\
 &&~~~~~~~~~~~~~~=|\dot f(y+(\bar \theta-\theta_*)\odot x) (\alpha-\alpha'+(\beta-\beta')x)|\\
 &&~~~~~~~~~~~~~~\leq \sup_{\theta\in \Phi}| \dot f(y+(\theta-\theta_*)\odot x)| (|\alpha-\alpha'|+|\beta-\beta' ||x|).
 \end{eqnarray*}
 From condition (R) ii) there thus exists a nonnegative constant $c$ such that
  \begin{eqnarray*}
 &&\int_\R  |\Psi_\theta(y)-\Psi_{\theta'}(y)|\\
 &&\leq \sum_{j=0,1} (|\alpha-\alpha'|^{1-j}+| \beta-\beta'|^{j})\\
 &&\int_{\R\times \R}|x|^ j(\sup_{\theta \in \Phi}|\dot f(z+(\theta-\theta_*)\odot x)|+\sup_{\theta \in \Phi}|\dot f_0(z+\theta\odot x)|) h(x) dx dz\\
&& \leq c\|\theta-\theta'\|_2.
\end{eqnarray*}
\subsection{\it Uniform  Lipschitz property of  $H(\cdot,\vartheta)$}\label{integ_lipH}
Let us write
\begin{eqnarray*}
&&H(y,\vartheta)-H(y,\vartheta')\\
&=&\frac{1}{p}\left(F_\theta(y)-F_{\theta'}(y)+F_\theta(-y)
-F_{\theta'}(-y)\right)+\frac{p-p'}{pp'}(F_{\theta'}(y)+F_{\theta'}(y))\\
&&+\frac{1-p}{p}\left(J_\theta(y)-J_{\theta'}(y)+J_\theta(-y)-J_{\theta'}(-y)\right)+\frac{p-p'}{pp'}(J_{\theta'}(y)+J_{\theta'}(-y)).
\end{eqnarray*}
To prove the uniform Lipschitz property of $H(\cdot,\vartheta)$ we need to prove it for  $J_\theta(\cdot)$ and $F_\theta(\cdot)$.  We begin with the simplest
term $J_\theta(\cdot)$.  According again to the mean value theorem, for all $y\in \R$,  all $(x,z)\in \R^2$ with $z\leq y$, and  all  $(\theta,\theta')\in  \Phi^2$ there exists $\bar \theta:=\bar \theta(x,z,\theta,\theta')$ belonging to the line segment with extremities
$\theta$ and $\theta'$ such that
\begin{eqnarray*}
|J_\theta(y)-J_{\theta'}(y)|&\leq& \int_{-\infty}^y \int_\R |f_0(z+\theta\odot x) -f_0(z+(\theta'\odot x) |h(x)dxdz\\
&=& \int_{-\infty}^y \int_\R |\dot f_0(z+\bar \theta\odot x)(|\alpha-\alpha'|+|\theta-\theta' ||x|)h(x)dxdz,\quad\quad  \\
&\leq&\int_{-\infty}^y \int_\R \sup_{\theta \in  \Phi} |\dot f_0(z+\theta\odot x)h(x)dxdz|\alpha-\alpha' |\\
&&  +\int_{-\infty}^y \int_\R |x|\sup_{\theta \in  \Phi} |\dot f_0(z+ux)h(x)dxdz|\beta-\beta' |\\
&\leq& c \|\theta-\theta' \|_2,
\end{eqnarray*}
where $c$ denotes a nonnegative constant arising from condition (R) ii). Using the same kind of argument we prove that there exists a nonnegative contant $ c'$ such that for all $(y,\theta,\theta')\in \R\times \Phi^2$
$$|F_\theta(y)-F_{\theta'}(y)|\leq c' \|\theta-\theta'\|_2.$$
In conclusion, for  all $y\in \R$, there exists a nonnegative constant $c''$ such that for all $(y,\vartheta,\vartheta')\in \R\times \Theta^2$
\begin{eqnarray*}
|H(y,\vartheta)-H(y,\vartheta')|\leq \frac{2}{\delta}(c+c')\|\theta-\theta' \|_2+4\frac{|p-p'|}{\delta^2}
\leq c'' \|\vartheta-\vartheta\|_3.
\end{eqnarray*}
\subsection{\it  Uniform almost sure rate of convergence of   $d_n$}\label{unif_rated}
Let us  consider
\begin{eqnarray*}
|d_n(\vartheta)-d(\vartheta)|\leq T_{1,n}(\vartheta)+T_{2,n}(\vartheta),
\end{eqnarray*}
where
\begin{eqnarray*}
T_{1,n}(\vartheta)&:=&\left |  \frac{1}{n}\sum_{i=1}^n H^2(V_i; \vartheta,\tilde F_{n,\theta}, \hat J_{n,\theta}) -H^2(V_i; \vartheta, F_{\theta},  J_{\theta})\right |,\\
T_{2,n}(\vartheta)&:=& \left | \frac{1}{n}\sum_{i=1}^n H^2(V_i; \vartheta,F_{\theta},  J_{\theta})-E\left(H^2(V_1; \vartheta,F_{\theta}, J_{\theta})\right) \right |.
\end{eqnarray*}

\noindent{\it  Uniform almost sure rate of convergence of   $T_{1,n}$.}
Note first that from boundedness of $H(\cdot; \vartheta,\tilde F_{n,\theta}, \hat J_{n,\theta})$ and $H(\cdot; \vartheta, F_{\theta},  J_{\theta})$ 
given by (\ref{boundH}), there exist  nonnegative constants $C$ and $C'$  such that 
\begin{eqnarray*}
T_{1,n}(\vartheta)&=& \left |  \frac{1}{n}\sum_{i=1}^n( H(V_i; \vartheta,\tilde F_{n,\theta}, \hat J_{n,\theta})+H(V_i; \vartheta, F_{\theta},  J_{\theta}))\right.\\
&&\left. \times (H(V_i; \vartheta,\tilde F_{n,\theta}, \hat J_{n,\theta})-H(V_i; \vartheta, F_{\theta},  J_{\theta}))\right |\\
&\leq& C\sup_{y\in \R}| H(y; \vartheta,\tilde F_{n,\theta}, \hat J_{n,\theta})-H(y; \vartheta, F_{\theta},  J_{\theta}))|\\
&\leq& C' \left(\sup_{y\in \R}| \hat J_{n,\theta}(y)-J_{\theta}(y)|+\sup_{y\in \R}|\tilde F_{n,\theta}(y)-F_{\theta}(y)|\right).
\end{eqnarray*}
Let us now denote
\begin{eqnarray*}
T^{(1)}_{1,n}&:=&\sup_{\theta\in \Phi}\sup_{y\in \R}| \hat J_{n,\theta}(y)-J_{\theta}(y)|,\\
T^{(2)}_{1,n}&:=&\sup_{\theta\in \Phi}\sup_{y\in \R}| \tilde F_{n,\theta}(y)-F_{\theta}(y)|.
\end{eqnarray*}

\noindent{\it Convergence rate of $T^{(1)}_{1,n}$.} 
For simplicity we will suppose that $\mbox{proj}_2(\Phi)\subset [0,A]$, where $A$ is a nonnegative real number and for all $(x,y)\in \R^2$, $\mbox{proj}_2: (x,y)\mapsto y$.
Let us  introduce $P^X_{n}=n^{-1} \sum_{i=1}^n \delta_{X_i}$  the empirical measure associated to the iid  sample $(X_1,\dots,X_n)$  with
common probability  distribution  $P^X$   with pdf  and cdf  resp. denoted  by  $h$  and $H$). We use the  functionnal notation  $Pf=\int f dP$. Notice now that,  according to expression (\ref{J-estim-final}),  we have  for all $y\in \R$:,
\begin{eqnarray}\label{emprocJ}
\hat J_{n,\theta}(y)-J_\theta(y)&=&\frac{1}{n}\sum_{i=1}^nF_0(y+\alpha+\beta X_i)-E(F_0(y+\alpha+\beta X))\\
&=&(P_n^X-P^X)F_0(y+\alpha+\beta ~\cdot).\nonumber
\end{eqnarray}
Let consider the class of functions
$$\mathcal{F}_0=\left\{  x\mapsto F_0(u+\beta x); \quad u\in \R ,\quad \beta\in [0,A]\right\}.$$
Since 
$$(P_n^X-P^X)F_0(u+\beta~\cdot)=(P_n^X-P^X)F_0(y+\beta (\cdot\vee 0))+(P_n^X-P^X)F_0(y+\beta (\cdot\wedge 0)),$$
it is enough to study the empirical process indexed by the classes of functions
\begin{eqnarray*}
\mathcal{F}^+_0&=&\left\{  x\mapsto F_0(u+\beta (x\vee 0)); \quad u\in \R ,\quad \beta\in [0,A]\right\},\\
\mathcal{F}^-_0&=&\left\{  x\mapsto F_0(u+\beta  (x\wedge 0)); \quad u\in \R ,\quad \beta\in [0,A]\right\}.
\end{eqnarray*} 
For simplicity we denote $\Gamma_{y,\alpha}(x)=F_0(y+\alpha (x\vee 0))$ and only consider the class $\mathcal{F}^+_0$, the class $\mathcal{F}^-_0$
being treated in a entirely same way. Since $F_0$ is a cdf, for $\beta_1\leq \beta\leq \beta_2$ and $u_1\leq  u\leq u_2$ we have
\begin{eqnarray*}
\Gamma_{u_1,\beta_1}(x)\leq \Gamma_{u,\beta}(x)\leq \Gamma_{u_2,\beta_2}(x),\quad\quad x\in \R,
\end{eqnarray*}
and, since  $F_0$ is supposed to be Lipschitz,
\begin{eqnarray*}
0 \leq  \Gamma_{u_2,\beta_2}(x)-\Gamma_{u_1,\beta_1}(x)\leq C(u_2-u_1+(\beta_2-\beta_1)(x\vee 0)). 
\end{eqnarray*}

Let consider now $\varepsilon>0$, and $(\overline{u_\varepsilon} , \underline {u_\varepsilon})\in \R^2$ such that
$$F_0(\overline{u_\varepsilon})\geq 1-\varepsilon,\quad \quad \mbox{and} \quad \quad F_0(\underline {u_\varepsilon})\leq \varepsilon.$$
Note that $\overline{u_\varepsilon}$ and $\underline {u_\varepsilon}$ do not depend on $\beta$. 
For all $N\in \N$, define
$$\underline {u_\varepsilon}=u_{1,\varepsilon}\leq u_{2,\v}\leq\ldots \leq u_{ N,\varepsilon}=\overline{u_\varepsilon},$$ and consider $N(\varepsilon)$ the smallest integer such that $u_{i,\v}-u_{i-1,\v}\leq \v$ for $i=2,\ldots,N(\v)$. We denote by $\lceil\cdot \rceil$ the integer part function. For all $\v$ small enough we clearly have
$$
N(\v)\leq \left\lceil\frac{\overline{u_\varepsilon}-\underline {u_\varepsilon}}{\v} \right\rceil\leq 2 \frac{\overline{u_\varepsilon}-\underline {u_\varepsilon}}{\v} .
$$
Let us now define $\alpha_{i,\v}=\v(i-1)$, $i=1,\ldots,M(\v)$, where $M(\v)=\left\lceil \lceil A+1 \rceil /\v  \right\rceil$ and thus   $\alpha_{M(\v),\v}>A$. Observe in addition that
\begin{eqnarray*}
\| \Gamma_{u_{i+1,\v}, \beta_{j+1,\v}}- \Gamma_{u_{i,\v}, \beta_{j,\v}}\|^2_{2,P^X}&=&c^2 E\left((u_{i+1,\v}- u_{i,\v} + (\beta_{j+1,\v}- \beta_{j,\v}) (X_1\wedge 0))^2\right)\\
&\leq& 2 c^2 \v^2+2C^2 \v E(X_1^2)\\
&=&2c^2\v^2\left(1+E(X_1^2)\right).
\end{eqnarray*}
Hence the expression  $$[ \Gamma_{u_{i+1,\v}, \beta_{j+1,\v}}- \Gamma_{u_{i,\v}, \beta_{j,\v}}  ] , \quad 1\leq i\leq N(\v), \quad 1\leq j\leq M(\v),$$ 
is a $\left(c\sqrt{2(1+E(X_1^2))}\right)$-covering of $\mathcal{F}_0^+$ in the $L_{2}(P^X)$-norm sense. Using the standard  notation  $N_{[]}(\cdot)$  (see  van der Vaart and Wellner \cite{VW96}) the covering number  of the class $\mathcal{F}_0^+$ is bounded as follows
\begin{eqnarray*}
N_{[]}\left(\v,\mathcal{F}_0^+,L_2(P^X)\right)\leq cN(\v)M(\v)\leq c'\frac{\overline{u_\varepsilon}-\underline {u_\varepsilon}}{\v^2}.
\end{eqnarray*}
Thus if there exist constants $C$ and $V$ such that 
\begin{eqnarray}\label{cond_bracket}
|\overline{u_\varepsilon} |\wedge |\underline {u_\varepsilon}|\leq C/\v^V,
\end{eqnarray}
 we get 
$N(\varepsilon)M(\v)\leq C/\v^{V+2}$ which  allows us to use 
Theorem 2.14.9,  p. 246  in  \cite{VW96} since their Condition (2.14.7), p. 245 is then satisfied   after replacing their constant $V$ by $V+2$.
Let us discuss condition (\ref{cond_bracket}). For $\v$ small enough this condition is true if  $\overline{y_\varepsilon}  \leq C/\v^V $ and 
$\underline {u_\varepsilon}\geq -C/\v^V$. Denoting by $F_0^{\leftarrow}$ the quantile function of $F_0$, condition (\ref{cond_bracket}) becomes
\begin{eqnarray}\label{cond_bracket_bis}
F_0^{\leftarrow}(1-\v)\leq C/\v^V\quad \quad \mbox{and}\quad \quad F_0^{\leftarrow}(\v) \geq -C/\v^V.
\end{eqnarray}
We  consider for simplicity  the first condition in (\ref{cond_bracket_bis}) (the second one being treated in the  same way);  it is equivalent to $F_0(C/\v^V)\geq 1-\v$, and taking $t=C/\v^V$ this condition turns into
\begin{eqnarray}
F_0(t)\geq 1-C/t^{1/V}.
\end{eqnarray}
Thus it suffices to have
\begin{eqnarray}
\liminf _{t\rightarrow \infty}\frac{-\log(1-F_0(t))}{\log(t)}>0.
\end{eqnarray}
Finally, using the symmetry of $f_0$,  condition (\ref{cond_bracket}) holds if 
\begin{eqnarray}\label{cond_bracket_fin}
\liminf _{t\rightarrow \infty}\frac{-2\log F_0(-t)}{\log(t)}>0,
\end{eqnarray}
which is insured by condition (R) vi).
In conclusion if (\ref{cond_bracket_fin}) is satisfied and $E(X_1^2)<\infty$ then, according to Theorem 2.14.16, p. 248 in van der Vaart and Wellner \cite{VW96},
we obtain
$$
\sup_{\theta\in \Phi }\|\ \hat J_{n,\theta}-J_\theta \|_\infty\leq \| P^X_n-P^X\|_{\mathcal F_0}=o_{a.s}(n^{-1/2+\gamma}),\quad \gamma>0.
$$

\noindent{\it Convergence rate of $T^{(2)}_{1,n}$}. 
Recall  that  $F_\theta $ is  the cdf of $Y_i-\theta\odot X_i$, {\it i.e.}
$$
F_\theta(y)=P(Y_i-\theta\odot X_i\leq y),
$$
and $$F_{n,\theta}(y)=\frac{1}{n}\sum_{i=1}^n\un_{Y_i-\theta \odot X_i\leq y}.$$
Let $K$ a kernel satifying (K). The $K$-regularized versions of $F_\theta$ and $F_{n,\theta}$  are 
$$\tilde F_\theta=K\ast F_\theta,\quad\quad \tilde F_{n,\theta}=K\ast F_{n,\theta}.$$
Let us denote by $P_n^{X,Y}$ the empirical measure 
$$
P^{X,Y}_n=\frac{1}{n} \sum_{i=1}^n \delta_{X_i,Y_i},
$$
and by $P^{X,Y}$ the law of $(X_1,Y_1)$.

The set of functions for which   $(x,y)\mapsto ax+by+c$  being a 3-dimensionnal vector space, Corollary 2.5 in Kuelbs and Dudley \cite{KD80} shows that
the class of sets 
$$
\mathcal{C}=\left\{ \left\{ (u,v)  \in \R^2:~ au+bv+c<0 \right\} ; ~~(a,b,c)\in \R^3\right\},
$$
is a Strassen log-log class, which implies that $a.s.$
\begin{eqnarray*}
\limsup_{n\rightarrow \infty} \sup_{\mathcal{C}} \sqrt{\frac{n}{2\log\log(n)}}(P_n^{Z}-P^{Z})(C)=\sup_{C\in \mathcal{C}}\sqrt{P^Z(C)(1-P^Z(C))}\leq 1/2.
\end{eqnarray*}
Since $\mathcal{C}$ contains the class
\begin{eqnarray*}
{\mathcal S}:=\left\{\left\{ (u,v)\in \R^2:~  v-(\alpha+\beta u)< y  \right\};~~(\alpha,\beta;y)\in \Phi\times\R\right\},
\end{eqnarray*}
it follows that, for all set $S\in{\mathcal S}$, $P^Z(S)=\int_{S} dP^{X,Y}(u,v)=P(Y-(\alpha+\beta X) <y)=F_{\theta}(y)$
and for the same reason $P_n^Z(S)=F_{n,\theta}(y)$,  we have 
\begin{eqnarray}\label{maj_loglog_fdr}
\limsup_{n\rightarrow \infty} \sup_{(\theta,y)\in \Phi\times \R}\sqrt{\frac{n}{\log\log(n)}}(F_{n,\theta}-F_{\theta})(y)\leq 1/2\quad \quad a.s.
\end{eqnarray}
Now if we replace $F_{n,\theta}$ by its regularized version $\tilde F_{n,\theta}$ the approximation is controlled as follows,
\begin{eqnarray}\label{repl_smooth_cdf}
&&\tilde F_{n,\theta}(y)-F_{n,\theta}(y)\nonumber\\
&&~~~~~~~~~~~~~=\tilde F_{n,\theta}(y)-E( \tilde F_{n,\theta}(y))+E(\tilde F_{n,\theta}(y))-F_\theta(y)+F_\theta(y)-F_{n,\theta}(y)\nonumber\\
&&~~~~~~~~~~~~~=\tilde F_{n,\theta}(y)-E( \tilde F_{n,\theta}(y))-[F_{n,\theta}(y)-E(F_{n,\theta}(y))]\nonumber\\
&&~~~~~~~~~~~~~~+E(\tilde F_{n,\theta}(y))-F_{\theta}(y),
\end{eqnarray}
recalling that $E(F_{n,\theta}(y))=F_\theta(y)$. The first   term on  the  right hand side of (\ref{repl_smooth_cdf}) satisfies
\begin{eqnarray*}\label{maj_K_TV}
\tilde F_{n,\theta}(y)-E(\tilde  F_{n,\theta}(y))&=&\int_\R K(y-u)d (F_{n,\theta}-E(F_{n,\theta}))(u)\nonumber\\
&=&\int_\R (F_{n,\theta}-E(F_{n,\theta}))(u) dK(y-u)\nonumber\\
&=&\int_\R(F_{n,\theta}-E(F_{n,\theta}))(y-s)dK(s).\nonumber
\end{eqnarray*}
Thus, if we denote $\Delta_{n,\theta}(y):=F_{n,\theta}(y)-E(F_{n,\theta}(y))=F_{n,\theta}(y)-F_{\theta}(y)$, we obtain
\begin{eqnarray*}
&&\left| \tilde F_{n,\theta}(y)-E( \tilde F_{n,\theta}(y))-[F_{n,\theta}(y)-E(F_{n,\theta}(y))]\right|\\
&&~~~~~~~~~~~~~~~~~~~~\leq \left| \int_\R (\Delta_{n,\theta}(y-s)-\Delta_{n,\theta}(y))dK(s) \right|\nonumber\\
&&~~~~~~~~~~~~~~~~~~~~\leq \sup_{(\theta,y)\in \Phi\times \R}\left |   \Delta_{n,\theta}(y) \right| \left\| K\right\|_{TV}.
\end{eqnarray*}
The last bias-term on the right hand side of (\ref{repl_smooth_cdf}) can be studied using the $R_{2n}$ bound in \cite{SW86}, p. 766,  equation  (e),
which establishes that for each $\theta\in \Phi$
\begin{eqnarray}\label{maj_R2n}
\sup_{y\in \R} |E(\tilde F_{n,\theta})-F_\theta |(y)\leq \frac{\|\dot f_{\theta}\|_\infty k_2} {2}.
\end{eqnarray}
If $K$ is replaced by $K_n(\cdot)=K(\cdot/b_n)$ and we let  $\displaystyle k_{2,n}:=\int_\R x^2 dK_n(x)=b_n^2 k_2$, then (\ref{maj_loglog_fdr}--\ref{maj_R2n}) lead to 
\begin{eqnarray}\label{maj_loglog_global}
\limsup_{n\rightarrow \infty} \sqrt{\frac{n}{\log\log(n)}}\sup_{(\theta,y)\in \Phi\times \R}(F_{n,\theta}-F_{\theta})(y)<\infty \quad \quad a.s.,
\end{eqnarray}
whenever $\limsup (n/\log\log(n))^{1/2}k_{2,n}<\infty$ which holds when 
\begin{eqnarray}
\limsup_{n\rightarrow \infty} \sqrt{\frac{n}{\log\log(n)}}b_n^2<\infty,
\end{eqnarray}
and $\sup_{\theta\in \Phi}\|\dot f_{\theta}\|_\infty <\infty$ which has been proved in Section 5.3 under Condition (R) ii).\\

\noindent{\it  Uniform almost sure rate of convergence of   $T_{2,n}$.}
Considering for all $i\geq 0$, the random variable $W_i(\vartheta):=H^2(V_i;\vartheta)$, where $\vartheta\in \Theta$, we see that 
$$
\sup_{\vartheta\in \Theta} T_{2,n}(\vartheta)=\sup_{\vartheta\in \Theta}\left|\frac{1}{n}\sum_{i=1}^n W_i(\vartheta)-E(W_i(\vartheta))  \right|,
$$
where the right hand  term is the supremum of an empirical process indexed by a class of Lipschitz bounded functions, which is known to be $o_{a.s.}(n^{-1/2+\gamma})$
for all $\gamma>0$ (see \cite{BMV06}, for details), which concludes the proof.\\

\noindent {\bf Acknowledgments}. The author  thanks the referees for their helpful and constructive comments. He is also very gratefull to Philippe Barbe and David Hunter  for their help and good advice during the writing of this manuscript.

\end{document}